\documentclass[11pt]{article}

\usepackage[utf8]{inputenc}
\DeclareUnicodeCharacter{202F}{\,}

\usepackage{a4wide}
\usepackage{geometry} % (optional; keep if you plan to control margins here)

\usepackage{amsmath,amssymb,amsfonts,amsthm}
\numberwithin{equation}{section}

\usepackage{graphicx}
\usepackage{float}
\restylefloat{figure}
\usepackage{wrapfig}

\usepackage[table]{xcolor}
\usepackage{array,tabularx,longtable,booktabs,multirow,makecell}

\newcolumntype{s}{>{\columncolor{blue!10}}c}
\arrayrulecolor{gray!10}

\usepackage[caption=false]{subfig}
\usepackage{caption}
\usepackage{enumerate}
\usepackage{enumitem}

\usepackage{authblk}
\usepackage{titlesec}

\usepackage{varioref}
\usepackage{hyperref}
\hypersetup{
	colorlinks=true,
	linkcolor=violet,
	citecolor=cyan,
	urlcolor=teal
}

\def\RR{\hbox{I\kern-.2em\hbox{R}}}
\newcommand{\myqed}{\hbox to 0pt{}\hfill$\rlap{$\sqcap$}\sqcup$\vspace{3mm}}

\graphicspath{{Figure/}}
\date{}

\begin{document}
	\title{Climate-Driven Dengue Forecasting in Bangladesh: Division-Specific Feature-Set Design and Lag Structure}

	\author[1, 2]{\small Faizunnesa Khondaker\thanks{Email: kh.faizunnesa@gmail.com}}
	\author[1]{\small Md. Kamrujjaman\thanks{Corresponding author   Email: kamrujjaman@du.ac.bd} }
	%	\author[2]{\small  Md. Shahidul Islam\thanks{Email: mshahid@du.ac.bd}}

	\affil[1]{\footnotesize Department of Mathematics, University of Dhaka, Dhaka 1000, Bangladesh}
	\affil[2]{\footnotesize Department of Mathematics, Jagannath University, Dhaka 1100, Bangladesh}
%	\affil[3]{\footnotesize Department of Applied Mathematics and Data Science, Asian University for Women, Chattogram 4000, Bangladesh}
	\affil[ ]{\textbf{Corresponding author:} Md. Kamrujjaman, \href{mailto:kamrujjaman@du.ac.bd}{kamrujjaman@du.ac.bd}}
	\maketitle
	
	\vspace{-1.0cm}
	\noindent\rule{6.35in}{0.02in}\\
	\section*{Abstract}
			Bangladesh exhibits marked year-to-year variability in dengue, partly driven by meteorological fluctuations that shape \textit{Aedes} breeding-site persistence, mosquito development, and transmission. We exploit a contrast between Dhaka (consistently high burden) and Barishal (recently rising burden despite lower population density) and frame feature-set design and predictor structure as the main methodological contributions. Using monthly dengue data from DGHS \cite{DGHS} and meteorological data from World Weather Online \cite{Weather} for January 2022--October 2025, we compare four climate feature sets that vary wetness (rainy days vs.\ rainfall) and sunshine (sun days vs.\ sun hours), while temperature and humidity appear in all sets. We evaluate two predictor configurations: lagged climate covariates only, and lagged climate covariates plus 1-month lagged dengue incidence ($Y_{t-1}$). Climate lags (0--4 months) are applied in correlation and forecasting. Both divisions show similar delayed associations: rainfall metrics peak positively near a 2-month lag, humidity near a 1-month lag, sunshine metrics are most negative around a 2-month lag, and temperature is weakly positive at longer lags. We then benchmark MPR, ANN, XGBoost, and SARIMAX across all sets. Best performance differs: Dhaka favors ANN-1 with SET-1 (RMSE=2176.70, MAE=1282.00, MAPE=31.54\%), whereas Barishal favors SARIMAX(0,1,1)(1,0,0,12) with SET-2 (RMSE=817.56, MAE=717.78, MAPE=39.96\%). Analyses use consistent monthly aggregation and division-specific tuning.

	\noindent\textbf{Keywords:} Dengue; Time-series forecasting; SARIMAX; ANN; XGBoost; Poisson regression.\par
	\noindent\textbf{MSC 2020:} 92D30; 62M10; 62J12.\par
	
	\noindent\rule{\textwidth}{0.4pt} % optional separator

	\section{Introduction}
	Dengue is a major mosquito-borne viral disease that has expanded beyond its traditional tropical-subtropical range with globalization and climate change \cite{Lee2013,Simo2019}. It is caused by four antigenically distinct serotypes (DENV-1--DENV-4) and is transmitted mainly by the urban vector \textit{Aedes aegypti} and secondarily by \textit{A. albopictus} \cite{Ong2021}. Now endemic in nearly 100 countries, dengue places over half of the world's population at risk and causes an estimated 50-200 million infections annually, with severe dengue contributing substantial mortality in endemic settings \cite{Murray2013}. Transmission is strongly shaped by climate and human environments because \textit{Aedes} mosquitoes breed in small, often human-made water containers and pass through aquatic (egg-larvae-pupae) and adult stages \cite{Life} that respond to temperature, rainfall, and humidity via breeding-site availability, development speed, survival, and biting activity \cite{Nasirian2025}. Higher temperature and humidity favor faster mosquito development and more biting \cite{Turner2015}, while rainfall creates standing water that supports \textit{Aedes} breeding. Sunny days can contribute indirectly by altering evaporation and the persistence of small water habitats after rainfall, while dry spells may increase water storage that creates additional container breeding sites. In Bangladesh, monsoon rainfall variability, warm humid seasons, dense urbanization, and water-storage practices together create favorable conditions for \textit{Aedes} proliferation and dengue spread.
	
	Dengue-climate associations have been investigated using time-series methods \cite{Islam2021}, generalized linear models \cite{Hossain2022}, Poisson regression \cite{Hossain2023,Karim2023}, negative binomial regression \cite{Muurlink2018}, fuzzy analysis \cite{Shova25}, machine-learning approaches \cite{Dey2022}, and data-driven \cite{mobin24}. \cite{Alam2025} shows that climate strongly shapes dengue dynamics in Bangladesh and calls for meteorology-based early warning with adaptive, climate-informed control and surveillance. \cite{Chowdhury2025} combines climate-driven trend analysis with ANN and XGBoost forecasting to identify the best predictive model and key climatic drivers for early warning and control. In terms of model comparisons, \cite{Akter2024} reports that a neural network model outperforms ARIMA-based time series analysis for predicting dengue outcomes in Dhaka, achieving markedly lower error (e.g., MAPE) than the best ARIMA specification. In \cite{Hossain2024} economic time-series methods (ADF, VAR/VECM, Granger causalty, IRF, and variance decomposition) to quantify meteorological effects, finding rainfall as the dominant driver with a unidirectional causal link to dengue and both short and long-run impacts, supporting climate-based early warning and strengthened monitoring in endemic urban areas. In parallel, artificial neural networks (ANNs) are well suited for capturing nonlinear relationships and have shown strong performance in diverse applications \cite{Azeem2023,Bukhari2021}, while eXtreme Gradient Boosting (XGBoost) is also effective for nonlinear data and has achieved high predictive accuracy across multiple domains \cite{Noorunnahar2023,Guan2023}.\\
	\begin{figure}[H]
		\centering
		\includegraphics[width=5.5 in]{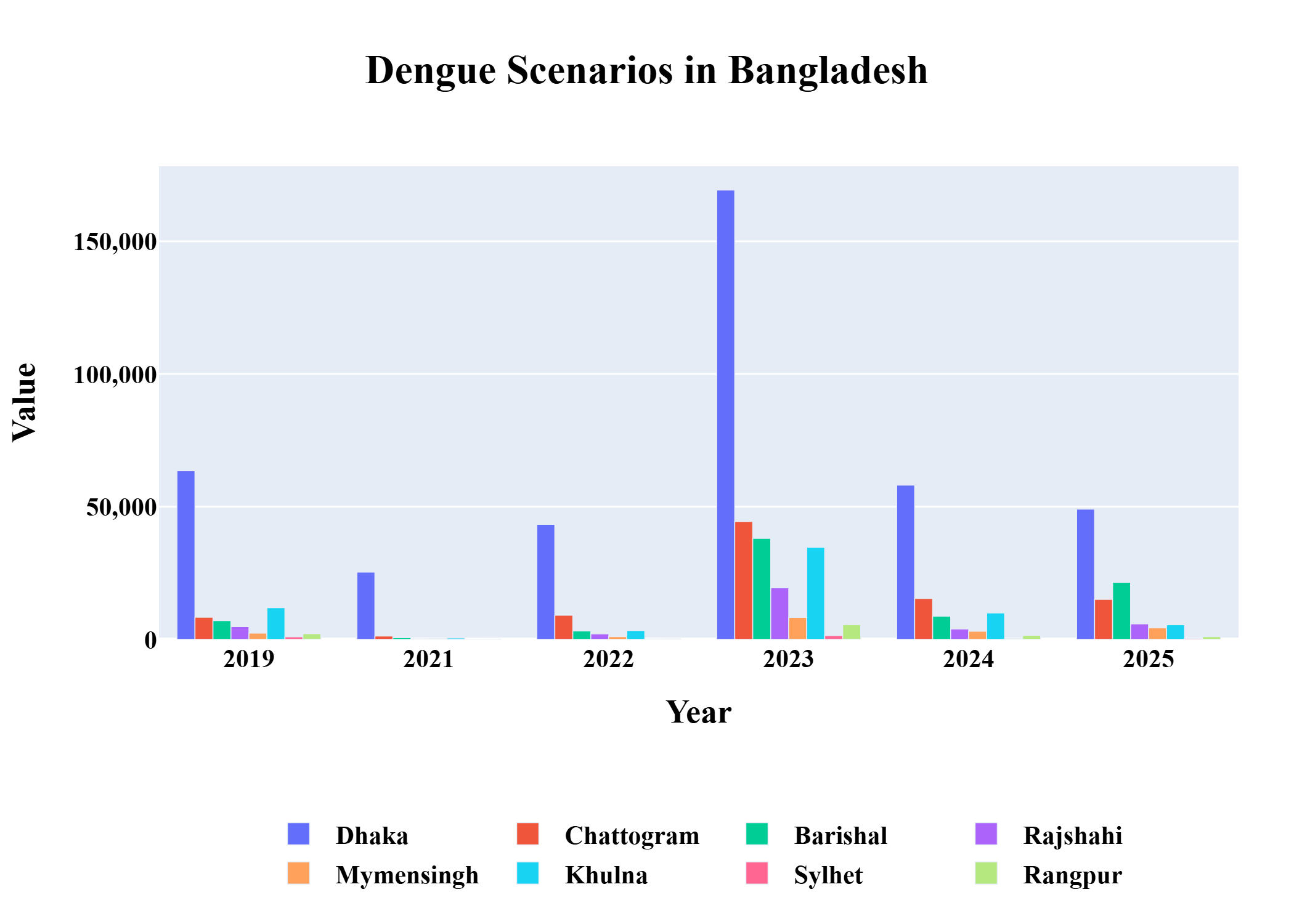}
		\caption{Division-wise annual dengue burden in Bangladesh (2019-2025), highlighting spatial heterogeneity and year-to-year variation across divisions.}	
		\label{div}
	\end{figure}
	Dengue transmission responds strongly to meteorological variability, because changes in seasonal conditions can alter how long breeding sites persist and how efficiently \textit{Aedes} mosquitoes develop and transmit infection; therefore, the burden may fluctuate markedly from year to year when weather patterns and local conditions coincide. 
	\begin{figure}[H]
		\centering
		\includegraphics[width=4 in]{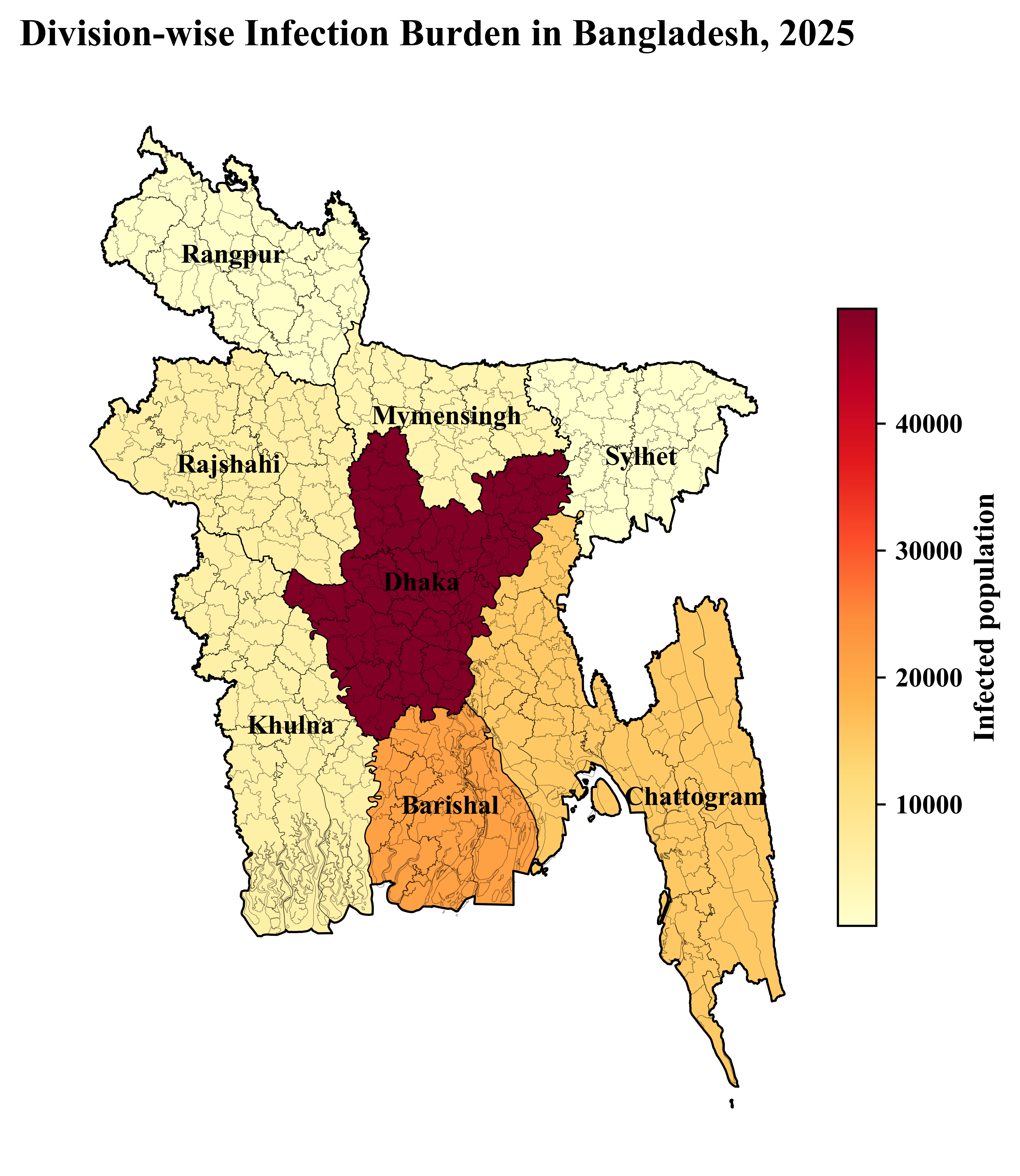}
		\caption{Division-wise heatmap of infected population in Bangladesh in 2025, illustrating the spatial variation in infection burden across divisions, with darker shades indicating higher case counts.}	
		\label{heat}
	\end{figure}
	From 2019 to 2025, the dengue picture across the eight divisions of Bangladesh--Barishal, Chattogram, Dhaka, Khulna, Mymensingh, Rajshahi, Rangpur, and Sylhet-- reads like a shifting map of risk (Figure~\ref{div}). After the high counts in 2019, the numbers fall sharply in 2021, then rebound in 2022, before the country experiences a dramatic nationwide surge in 2023. The following year (2024) shows a clear retreat from that peak in every division. But by 2025, the pattern is no longer a uniform decline--most divisions continue to come down, while one division moves strongly in the opposite direction \cite{DGHS}. That divergence is most striking in Barishal: it jumps from 8,763 cases in 2024 to 21,547 in 2025, even though Barishal is the least populous division, with a population of about 9.10 million and it also has the lowest population density at about 688 people per km$^2$. At the same time, Dhaka continues to dominate the national burden: it records 49,070 cases in 2025 (the highest among all divisions) (Figure~\ref{heat}). It also continues to attract inward migration and concentrates population within $20,508.8$ km$^2$, accommodating about $44.18$ million residents and the country's highest density of about $2,156$ people per km$^2$ \cite{BBS}. Taken together, Dhaka and Barishal capture the central contrast motivating this study. Dhaka reflects a consistently high-burden setting where large population size and very high density can sustain transmission, whereas Barishal signals an emerging concern, showing a pronounced rise while most other divisions decline. 
	
	Building on this contrast, the remainder of this paper is organized as follows: Section ~\ref{Temporal} describes the temporal and seasonal distribution of dengue cases and the associated climate patterns. Section ~\ref{Framework} presents the data sources and outlines the overall modeling framework. Section ~\ref{Methodology} details the methodology, including data alignment and lag construction, decomposition, and correlation analyses, and the forecasting models. Finally, Section ~\ref{Discussion} summarizes the main findings, discusses their public-health implications, and provides concluding remarks.

	\section{Temporal and Seasonal Distribution of Dengue Cases}\label{Temporal}
\begin{figure}[H]
		\centering
		\subfloat[]{\includegraphics[width=4.2in,height=3.5in,keepaspectratio]{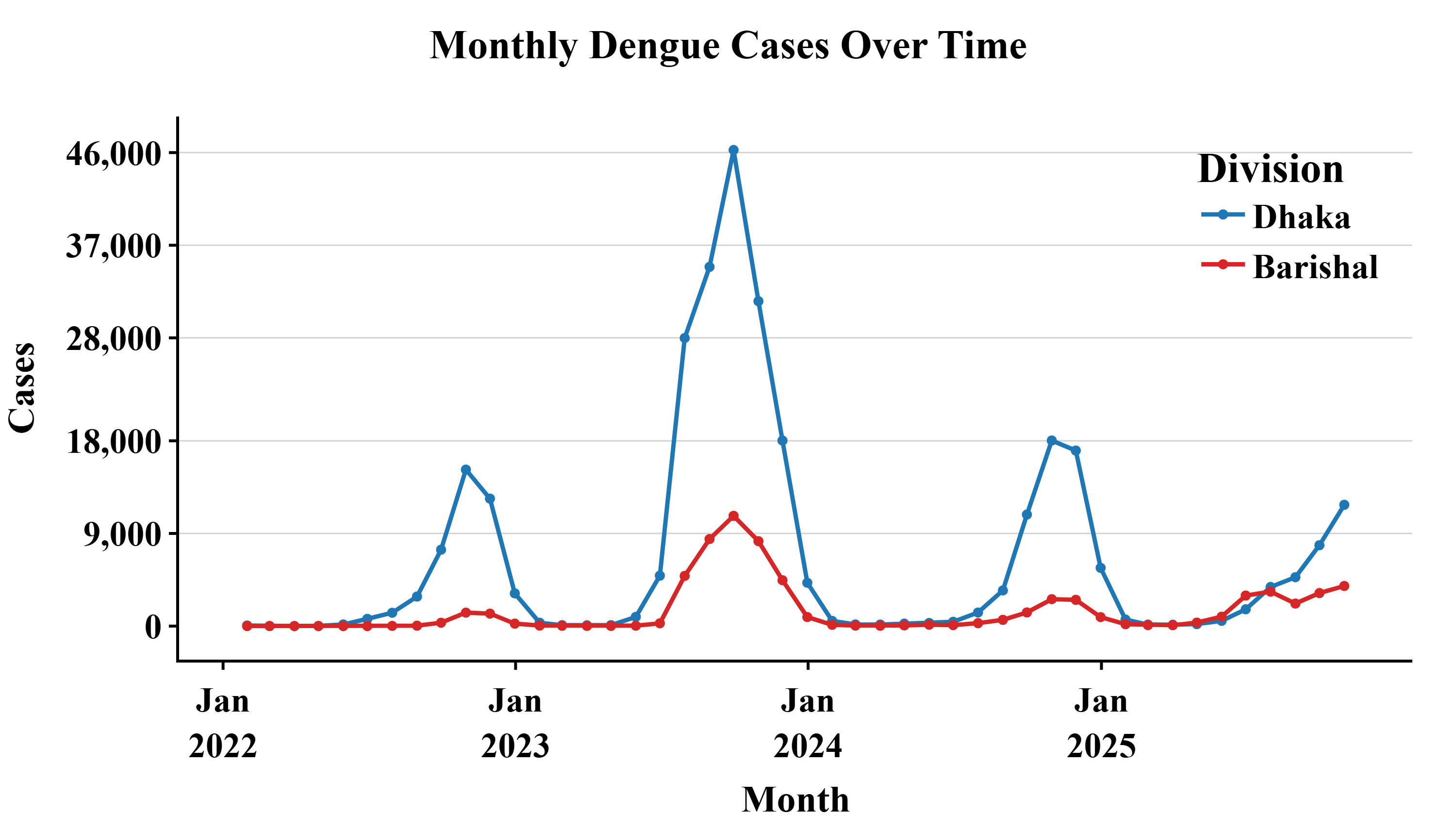}}\\
		\subfloat[]{\includegraphics[width=4.2in,height=3.5in,keepaspectratio]{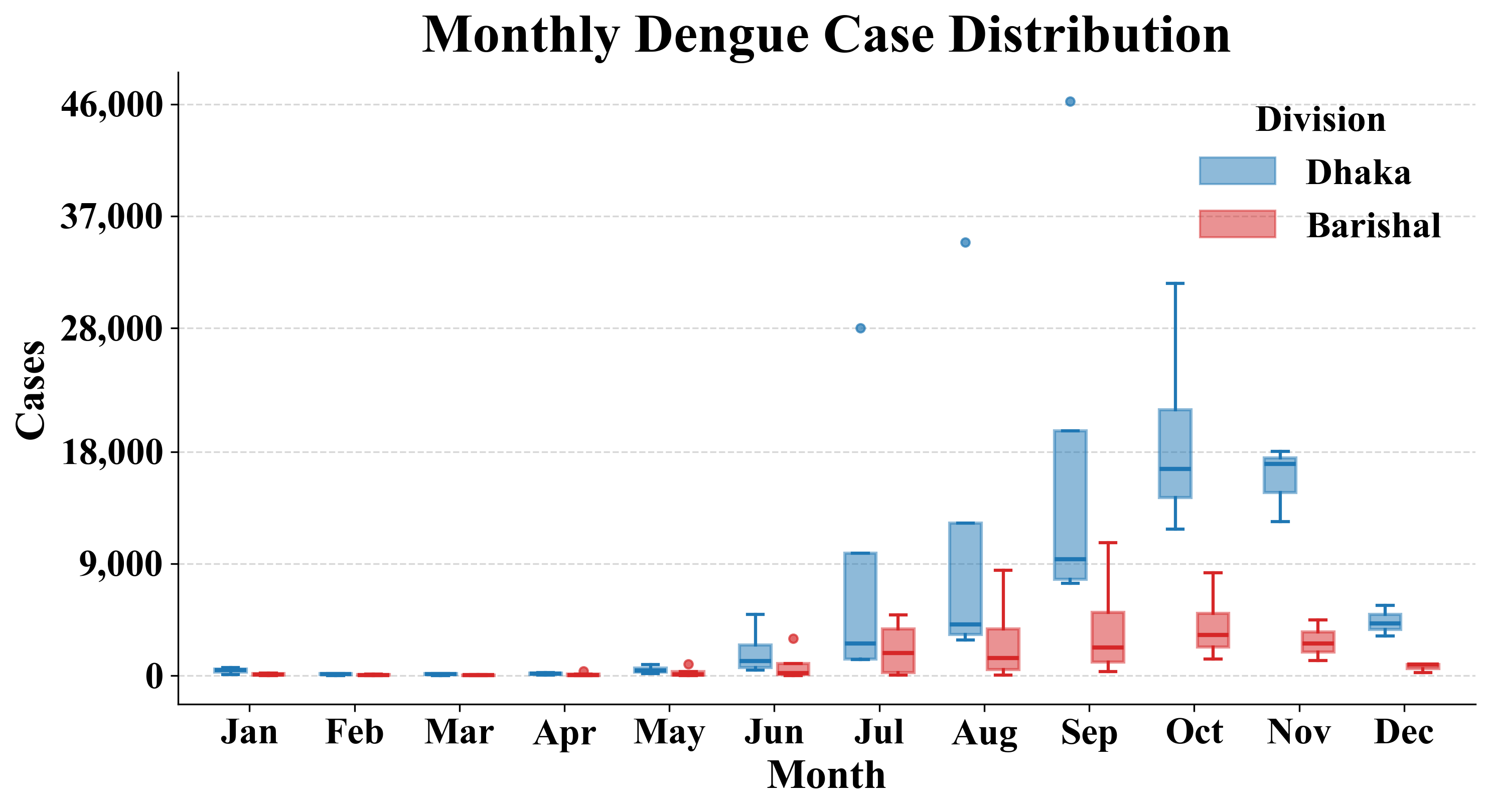}}
		\caption{Dengue cases in Dhaka and Barishal: (a) monthly reported cases from Jan 2022-- Oct 2025; (b) month-wise boxplots comparing seasonal dengue distributions between divisions.}	
		\label{case}
	\end{figure}

	Dengue transmission in Bangladesh is strongly seasonal, typically increasing during the monsoon and post-monsoon months when environmental conditions favor mosquito breeding. Consistent with this, dengue cases in Dhaka and Barishal division remain low in January-May and rise from June, peaking around September-October. As seen in Figure~\ref{case}(a), Dhaka consistently reports higher incidence, including a major late-2023 peak of about 46,250 cases with 10,712 cases on the same month September in Barishal, with smaller peaks in Dhaka during late 2022 and late 2024, while Barishal shows corresponding increases at a much lower level, and both divisions exhibit an upward trend again toward late 2025. Figure~\ref{case}(b) confirms the same seasonal pattern: Dhaka has much higher medians and wider spread (with extreme peak-season outliers), while Barishal follows the same timing but at a markedly lower magnitude \cite{DGHS}.

 	We consider six climate variables because they represent complementary mechanisms that influence dengue transmission, and their relative predictive contribution may vary across divisions. Rainy days measures the frequency of rainfall events, which can maintain persistent wetness and repeatedly refill container habitats, whereas rainfall amount measures the volume of water, which can increase standing water but may also flush larvae during heavy downpours. Sun hours reflects the intensity of sunshine, which promotes drying and can reduce local humidity, while sun days reflects the persistence of sunny conditions, indicating sustained dry periods. In addition, average temperature (Temp\_Avg) reflects the thermal environment that regulates mosquito development and viral replication rates, and humidity represents atmospheric moisture that affects adult mosquito survival and activity. 
 	\begin{figure}[H]
 		\centering
 		\includegraphics[width=5.5 in]{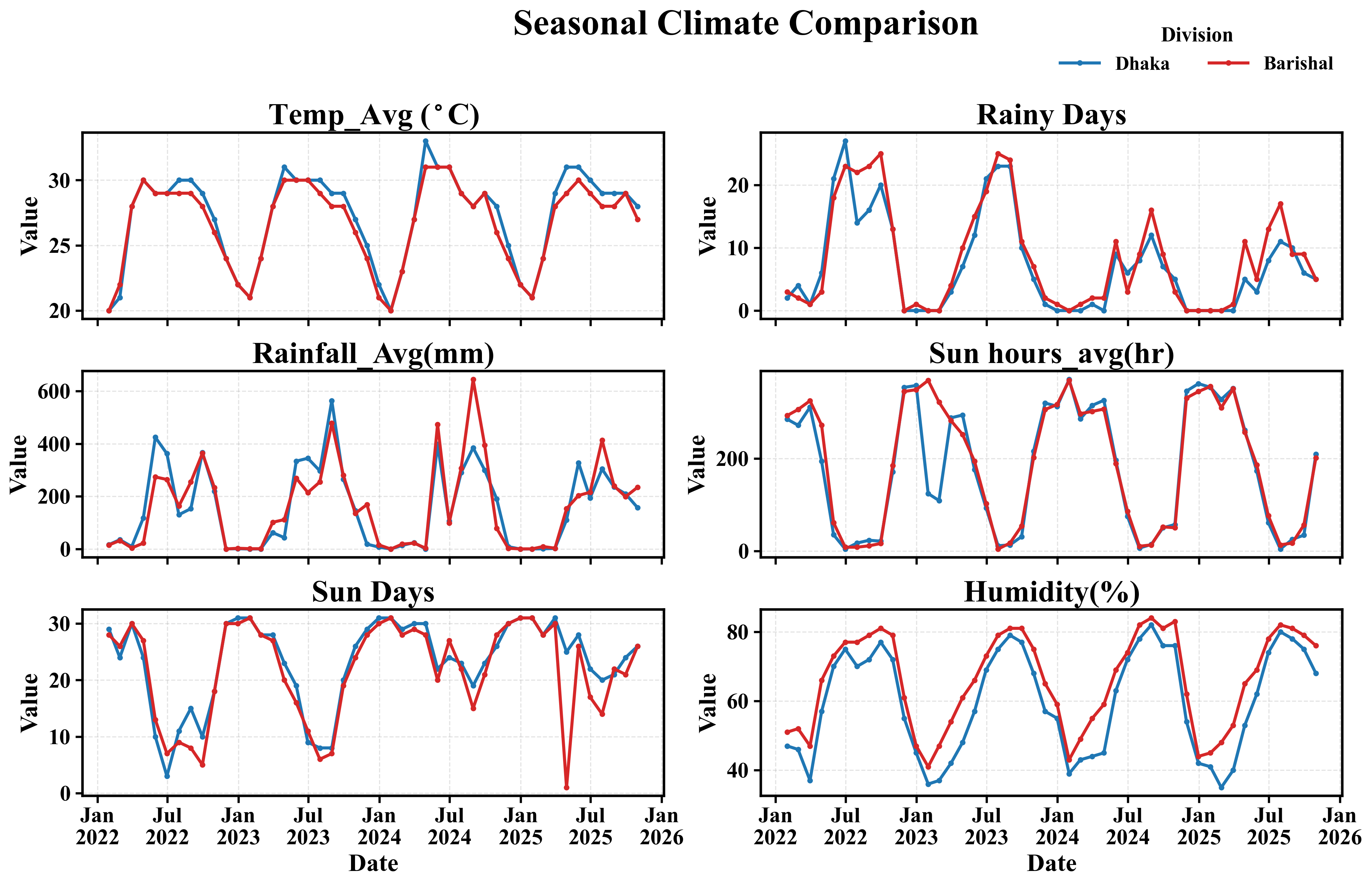}
 		
 		\caption{Time-series of key climate variables in the Dhaka and Barishal divisions over the study period.}	
 		\label{series}
 	\end{figure}
 \begin{figure}[H]
 	\centering
 \includegraphics[width=5.5 in]{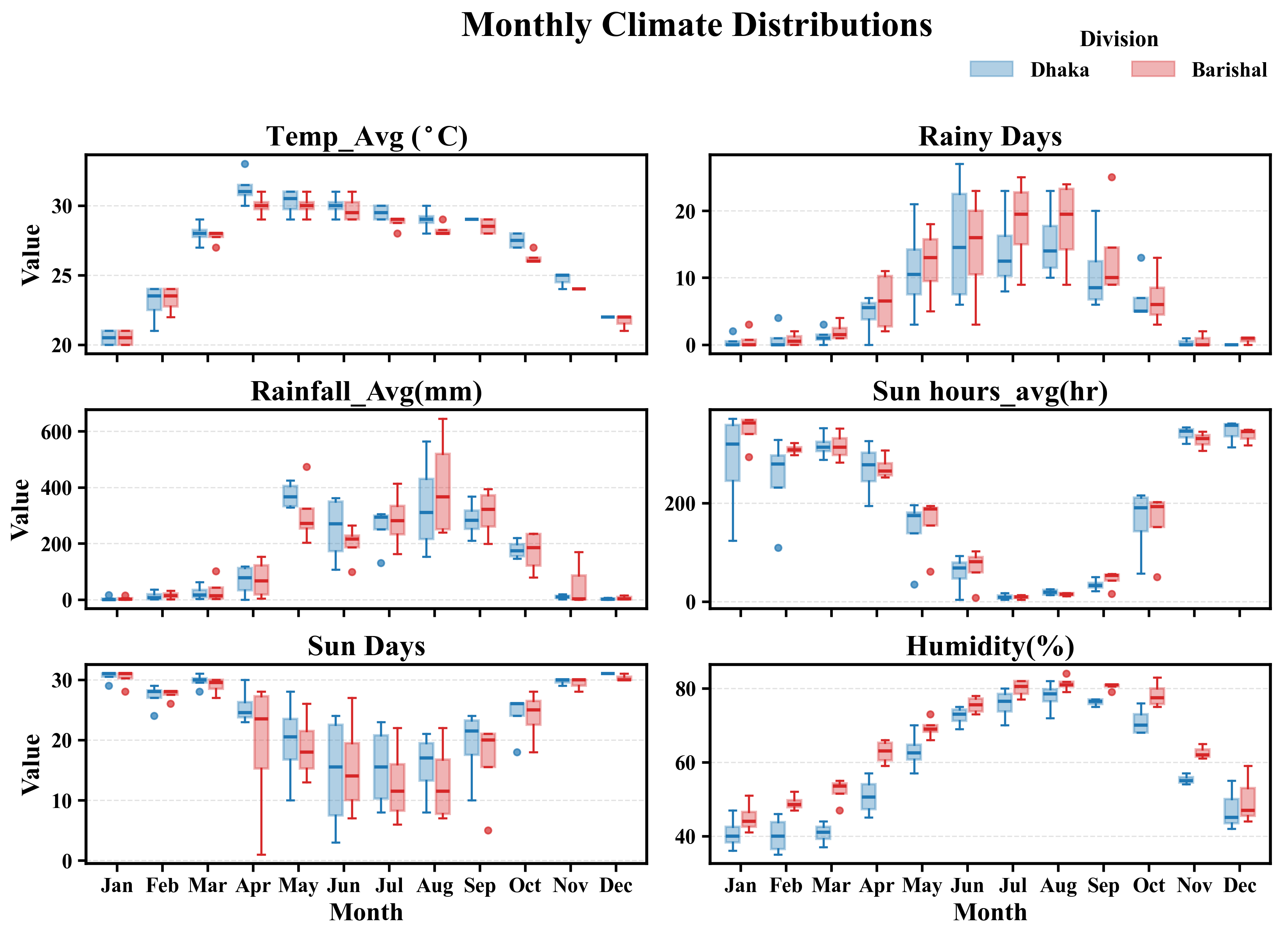}
 \caption{Month-wise boxplots of key climate variables in the Dhaka and Barishal divisions, summarizing monthly climate distributions across the study period.}	
 \label{boxtime}
\end{figure}
 	We observe from Figure~\ref{series} and ~\ref{boxtime} that Dhaka and Barishal experience a strongly synchronized, monsoon-driven climate cycle. Temperature increases from about 20--22$^\circ$C in winter (January) to roughly 30--33$^\circ$C during late spring-summer (April-August) and then declines toward year end, while humidity rises from approximately 40--55\% (January--March) to ~75--85\% in the monsoon peak (June--September). Over the same period, rainy days shift from near 0--2 days/month in the dry season to about 15--25 days/month in May-August, and rainfall increases from typically $<$50 mm/month in dry months to around ~150--350 mm/month in May--September, with occasional extremes reaching ~500--650 mm/month. In contrast, sunshine drops sharply during the monsoon: sun hours are generally high outside the monsoon (about ~250--350 hr/month) but fall to roughly ~0--80 hr/month in June--August, and sun days decrease from around ~28--31 days (dry months) to about ~10--20 days (monsoon months) \cite{Weather}. The time-series further show that the monsoon timing repeats each year (mid-year rainfall/humidity peaks and sunshine minima), whereas the boxplots compactly summarize the month-wise distributions and variability across the full period. Notably, the largest peaks observed in Dhaka and Barishal occur in September 2023, following the core monsoon period and coinciding with high rainfall/rainy days and humidity together with reduced sunshine, conditions that favor mosquito proliferation and typically lead dengue incidence be a short lag.
\\
To test alternative representations of wetness and sunshine, we define four sets: SET-1 (Temp + Rainy Days + Sun Hours + Humidity), SET-2 (Temp + Rainy Days + Sun Days + Humidity), SET-3 (Temp + Rainfall + Sun Days + Humidity), and SET-4 (Temp + Rainfall + Sun Hours + Humidity). Comparing results across SET-1 to SET-4 allows us to identify which wetness-sunshine combination provides the most informative exogenous inputs for dengue forecasting.
\section{Data Sources and Modeling Framework}\label{Framework}
	This study combines dengue surveillance records with meteorological observations for two administrative regions of Bangladesh: Dhaka division and Barishal division. The study period spans January, 2022 to October, 2025. Dengue case data were collected from the official daily dengue status reports published by the Directorate General of Health Services (DGHS), Bangladesh \cite{DGHS}, and meteorological data are obtained from WorldWeatherOnline \cite{Weather}. All analyses use secondary data: the datasets are already compiled and cleaned by the respective providers, and therefore the dowloaded datasets are treated as the final (raw) inputs for this study. No additional external processing of the original sources is required beyond aligning dates and constructing lagged variables for modeling.
	
	\subsection{Dengue case data}
	The outcome variable is the monthly dengue case count for each division (Dhaka and Barishal). DGHS \cite{DGHS} provides daily dengue status reports, which are aggregated to monthly totals to match the temporal resolution used in correlation screening and forecasting models. Monthly aggregation also reduces short-term reporting fluctuations and supports consistent comparison with monthly climate indicators.
	
	\subsection{Meteorological covariates}
	Monthly meteorological covariates are compiled separately for Dhaka and Barishal from World Weather Online \cite{Weather}. The selected variables represent key environmental conditions that influence mosquito activity and seasonal dengue transmission by capturing thermal conditions, rainfall intensity and wet-day frequency, sunshine exposure, and atmospheric moisture. Specifically, the climate dataset includes average temperature (Temp\_Avg), number of rainy days (Rainy Days), average daily sunshine hours (Sun hours\_avg(hr)), number of sunny days (Sun Days), and average relative humidity (Humidity(\%)).
	\subsection{Data alignment and lag construction}
	Dengue case counts and meteorological covariates are aligned to a common monthly time index spanning January 2022 to October 2025. To capture delayed climate effects, lagged versions of each climate variable (0-4 months) are generated and used in the correlation analysis and in the forecasting models, including Multivariate Poisson Regression (MPR), Artificial Neural Network (ANN), eXtreme Gradient Boosting (XGBoost), and Seasonal Autoregressive Integrated Moving Average with eXogenous regressors (SARIMAX). Observations with missing values arising solely from the lag construction are removed to ensure consistent training and testing datasets. This harmonized alignment and lag framework supports a direct comparison of climate-dengue associations between Dhaka and Barishal under a consistent modeling setup. 
	
	\subsection{Decomposition}
	Seasonal-trend decomposition using LOESS (STL) is commonly applied in environmental and public-health time-series studies \cite{Polwiang2020} to separate an observed series into interpretable components. STL decomposes the data into: (i) a trend component that captures long-term, low-frequency variation, (ii) a seasonal component that represents recurringg within-period patterns, and (iii) a remainder component that contains the variation not explained by the trend or seasonal structure. This approach is valued for its simplicity, robustness (particularly under the robust option), and clear visualization of temporal dynamics.
	
	Let $\mathcal{Y}_t$ denote the observed series and $\mathcal{T}_t, \mathcal{S}_t,$ and $\mathcal{R}_t$ denote the trend, seasonal, and remainder components, respectively. The additive decomposition is 
	\begin{align}
		\mathcal{Y}_t=\mathcal{T}_t + \mathcal{S}_t + \mathcal{R}_t.
		\end{align}
	Here, $\mathcal{Y}_t$ represents the monthly dengue case count, and t denotes time measured in months.
	\subsection{SARIMAX model}
	We use the Seasonal Autoregressive Integrated Moving Average with eXogenous regressors (SARIMAX) model to forecast monthly dengue cases while capturing both temporal dependence and climate forcing. SARIMAX extends SARIMA by incorporatingg external covariates as regressors, so dengue incidence is explained by seasonal autoregressive dynamics together with climate-driven variability.
	
	Let $y_t$ denote monthly dengue cases and $\mathbf{x}_t$ denote a vector of climate predictors (including lagged terms). The model is 
	\begin{equation}
		\Phi(B^{s})\, \phi(B)\,(1-B)^{d}(1-B^{s})^{D}y_t
		=\Theta(B^{s})\,\theta(B)\,\varepsilon_t + \boldsymbol{\beta}^{\top}\mathbf{x}_t,
	\end{equation}
where $B$ is the backshift operator, $s=12$ represents annual seasonality in monthly data, $(p, d, q)$ and $(P, D, Q)$ are the non-seasonal and seasonal orders, and $\varepsilon_t$ is a zero-mean error term.
To capture delayed climate effects, we use lagged exogenous predictors and evaluate SARIMAX models under the four predefined climate feature-set configurations (SET-1--SET-4). We evaluate a grid of SARIMAX specifications with non-seasonal orders $(p, d, q)$ and seasonal orders $(P, D, Q, 12)$ (with $s= 12$ for monthly seasonality), estimate parameters by maximum likelihood, and generate test-period forecasts using the corresponding exogenous inputs. Model selection relies on predictive accuracy (RMSE, MAE, MAPE), residual diagnostics via the Ljung-Box test (lag 12), and information criteria (AIC, BIC); for each feature set, we retain the model with the lowest test RMSE and report its AIC, BIC, and Ljung--Box $p$-value.
\subsection{Multivariate Poisson Regression (MPR)}
We model monthly dengue cases as count data using a Poisson generalized linear model. Let $Y_t$ denote the number of dengue cases in month $t$. We assume
\begin{equation}
	Y_t\sim\mathrm{Poisson}(\mu_t),
	\qquad
	\ln(\mu_t)=\beta_0+\sum_{i=1}^{k}\beta_i x_{it}+c\,Y_{t-1},
\end{equation}
where $\mu_t$ is the expected case count, $x_{it}$ are the (lagged) climate predictors, $\beta_0, \beta_i$ are  regression coefficients, and $c$ captures short-term persistence through the lagged dengue term. 

For each predefined featute set, we fit two specifications: MPR-1, which includes climate predictors only, 
\begin{equation}
	\ln(\mu_t)=\beta_0+\sum_{i=1}^{k}\beta_i x_{it},
\end{equation}
and MPR-2, which augments the climate predictors with the lagged dengue term, 
\begin{equation}
		\ln(\mu_t)=\beta_0+\sum_{i=1}^{k}\beta_i x_{it}+c\,Y_{t-1}.
\end{equation}
Parameters are estimated by maximum likelihood, and forecasts for the test period are generated using the corresponding covariates. %Predictive performance is evaluating using RMSE, MAE, and MAPE, and the best model is selected based on the lowest test RMSE.
\subsection{Artificial Neural Networks (ANN)}
The ANN forecasting models use a feed-forward multilayer perceptron (MLP) with an input layer, one or two hidden layers, and a single output node. The network learns nonlinear relationships between dengue incidence and predictors through backpropagation, using the tanh activation function and the Adam optimizer. To improve training stability and reduce over-fitting, we use regularization and learning-rate tuning, with a sufficiently large maximum number of training iterations. 

We implement two ANN specifications for each predefined feature set: ANN-1 uses lagged climate predictors only (four inputs), whereas ANN-2 adds a lagged dengue term (case lag) to capture short-term persistence (five inputs). The output of the network is the predicted number of monthly dengue cases. Hidden-layer sizes are chosen by hyperparameter tuning and result in either a single hidden layer (e.g., 16 neurons) or two hidden layers (e.g., 32 and 16 neurons), depending on the feature set.

\subsection{eXtreme Gradient Boosting (XGBoost)}
We apply eXtreme Gradient Boosting (XGBoost) as a tree-based ensemble learning approach to forecast monthly dengue cases using climate-driven predictors. XGBoost is an optimized implementation of gradient boosting that builds a strong predictor by sequentially adding regression trees, where each new tree is trained to reduce the errors of the previous ensemble. Its regularized formulation helps control model complexity and improves generalization, which is useful for nonlinear and interaction effects commonly observed in climate-dengue relationships.

Given observations $\{(x_i, y_i)\}_{i=1}^{n}$, the prediction at boosting iteration $t$ is 
\begin{equation}
	\hat{y}_i^{(t)}=\hat{y}_i^{(t-1)}+f_t(x_i),
\end{equation}
where $f_t(x_i)$ denotes the regression-tree function learned at boosting iteration $t$ and evaluated at the feature vector $x_i$. XGBoost minimizes the regularized objective
\begin{equation}
	\mathrm{Obj}^{(t)}=\sum_{i=1}^{n} l\!\left(y_i, \hat{y}_i^{(t-1)}+f_t(x_i)\right)+\Omega(f_t),
\end{equation}
in which $l(y_i, \hat{y}_i)$ is the loss function measuring the discrepency between the observed value $y_i$ and the prediction $\hat{y}_i$ and $\Omega(f_t)$ penalizes model complexity.

In our implementation, we use the XGBoost regressor with squared-error objective. For each predefined feature set, we consider two variants: XGB-1, which uses lagged climate predictors only, and XGB-2, which additionally includes a lagged dengue term (case lag) to capture short term persistence. Hyperparameters are tuned using grid-search cross-validation by minimizing prediction error on the training data.

For all modeling approaches (SARIMAX, MPR, ANN, and XGBoost), we evaluate predictive accuracy on an out-of sample test period using RMSE, MAE, and MAPE, and the best-performing specification is identified based on the lowest test RMSE.
\section{Methodology}\label{Methodology}
Seasonal-Trend decomposition using LOESS (STL) is used to separate the monthly dengue series into trend, seasonal, and residual components for Dhaka and Barishal. The decomposition reveals clear annual seasonality in both divisions, with Dhaka showing a larger seasonal amplitude, while Barishal displays a distinct seasonal swing that becomes more pronounced toward the later years. 
\begin{figure}[H]
	\centering
	\includegraphics[width=5 in]{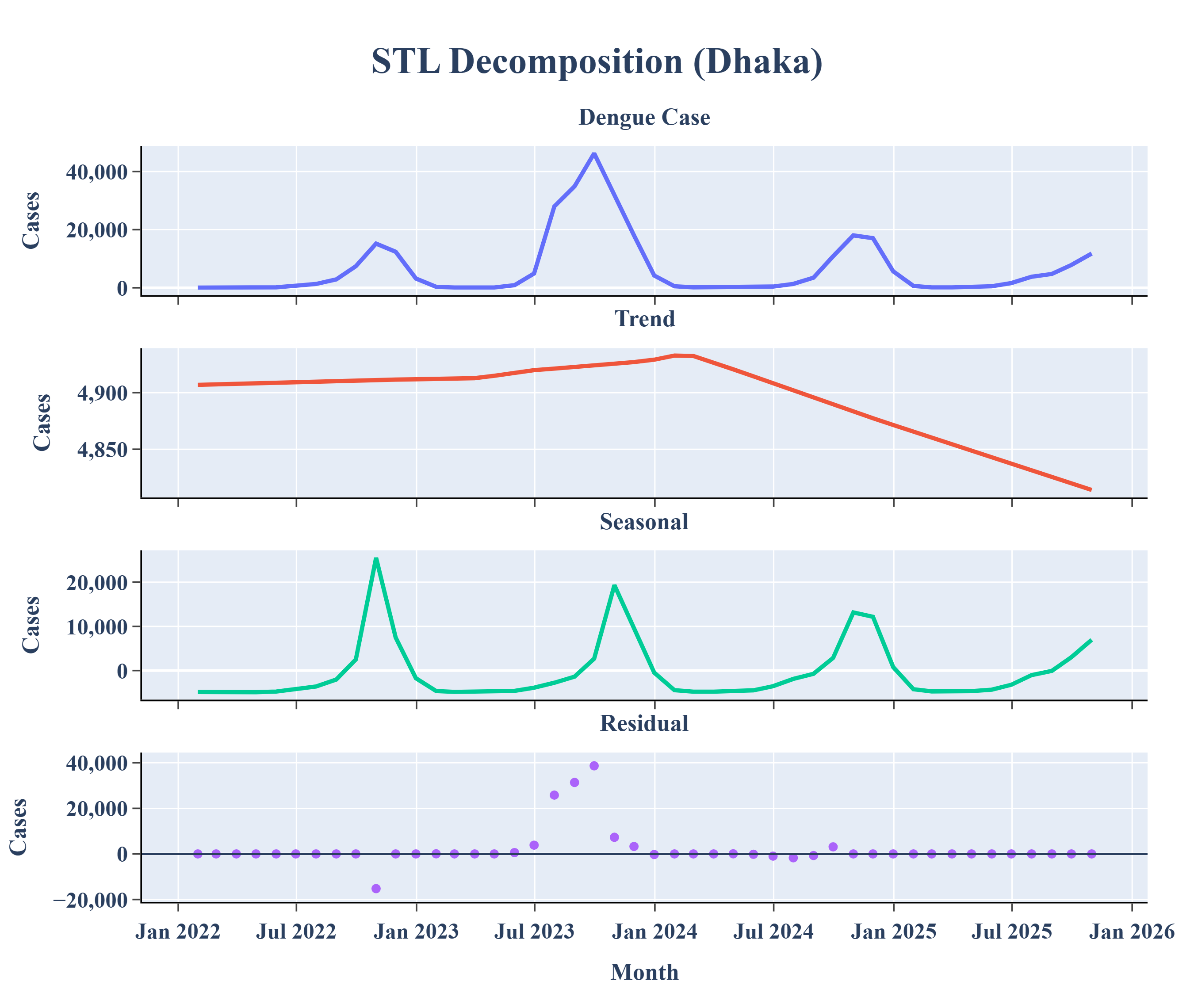}
	\caption{STL decomposition of monthly dengue cases for Dhaka, showing observed series, trend, seasonal component, and residuals.}	
	\label{stld}
\end{figure}
\begin{figure}[H]
	\centering
	\includegraphics[width=5 in]{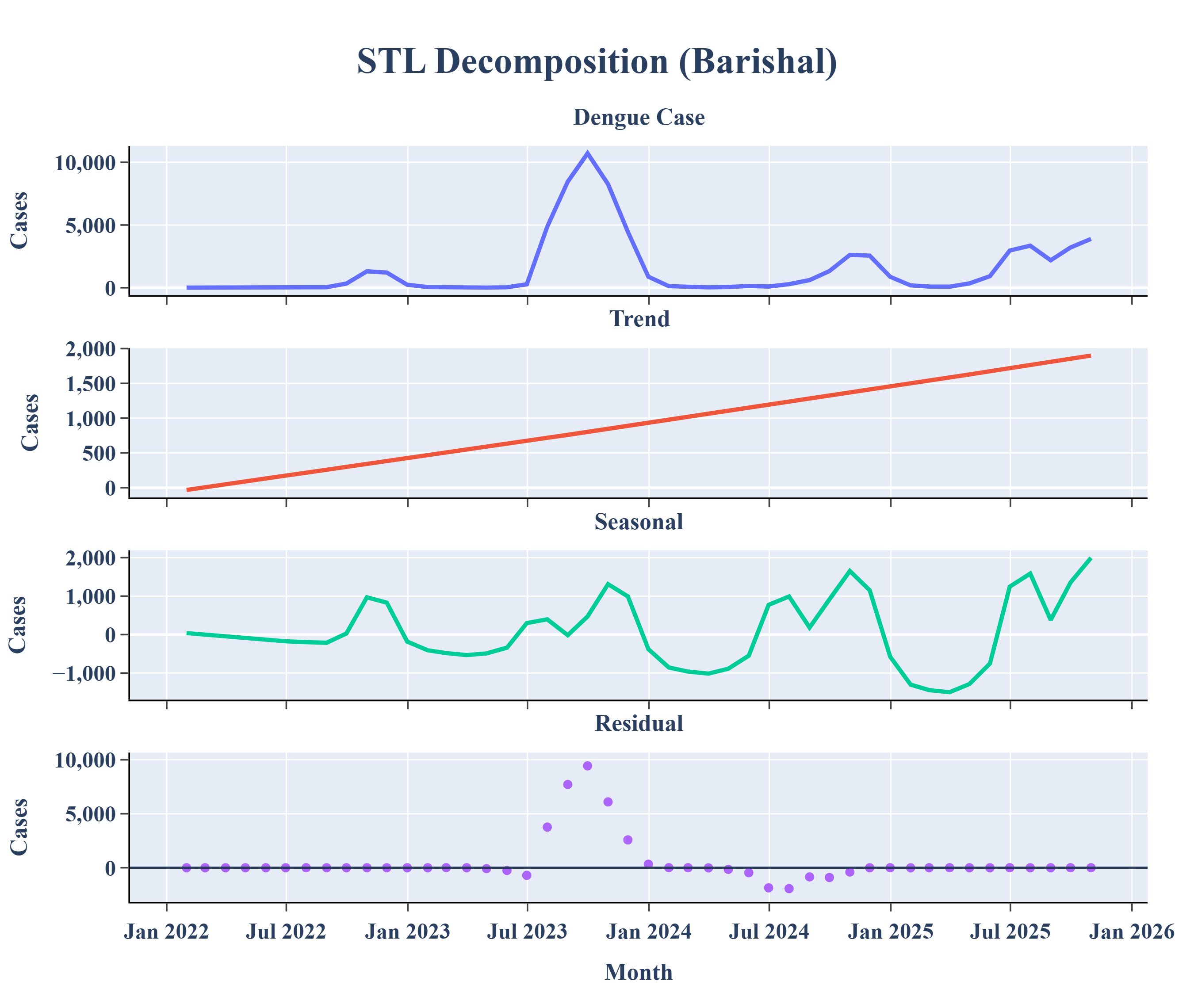}
	\caption{STL decomposition of monthly dengue cases for Barishal, showing observed series, trend, seasonal component, and residuals.}	
	\label{stlb}
\end{figure}
The trend in Dhaka increases slightly up to around 2023-early 2024 and then declines gradually (Figure~\ref{stld}), whereas Barishal exhibits a steadily increasing trend across the study period. Residuals are largest around the major 2023 outbreak in both divisions, indicating short-term shocks beyond the regular seasonal pattern (Figure~\ref{stlb}).
\begin{figure}[H]
	\centering
	\includegraphics[width=5 in]{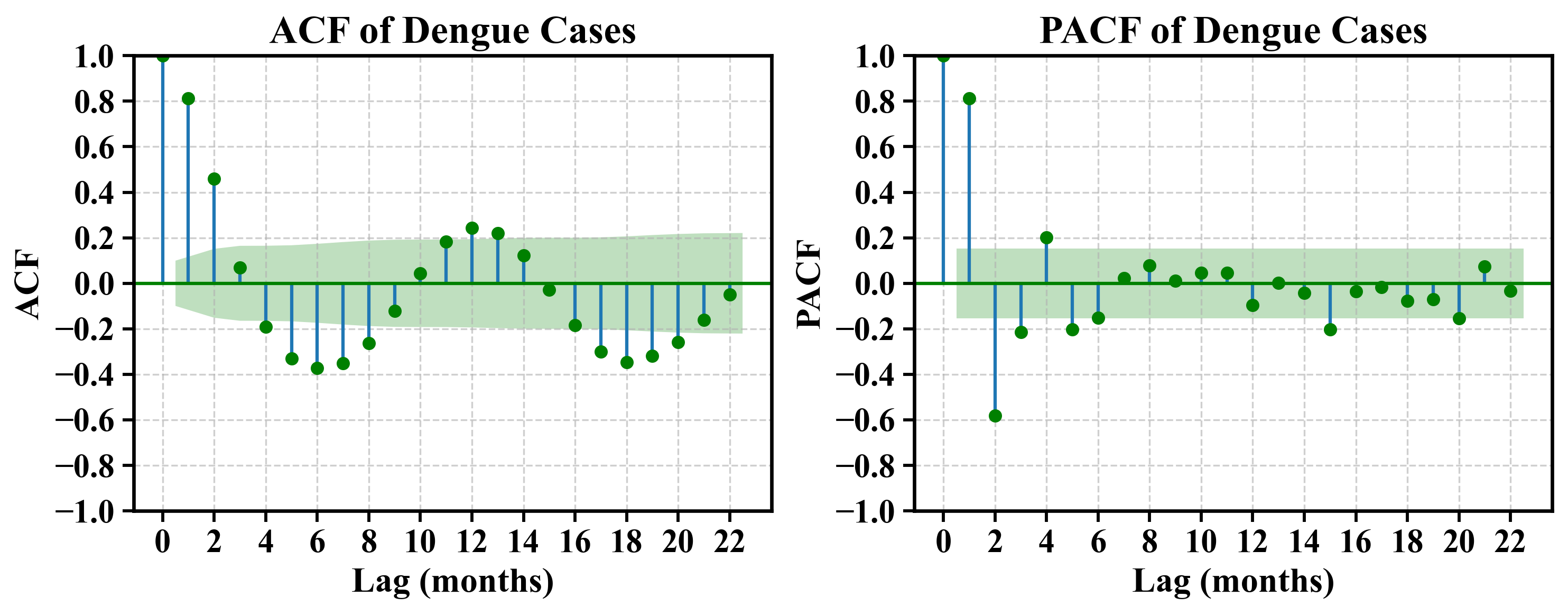}
	\caption{ACF and PACF plots of monthly dengue cases for Dhaka.}
	\label{acfd}
\end{figure}
\begin{figure}[H]
\centering
	\includegraphics[width=5 in]{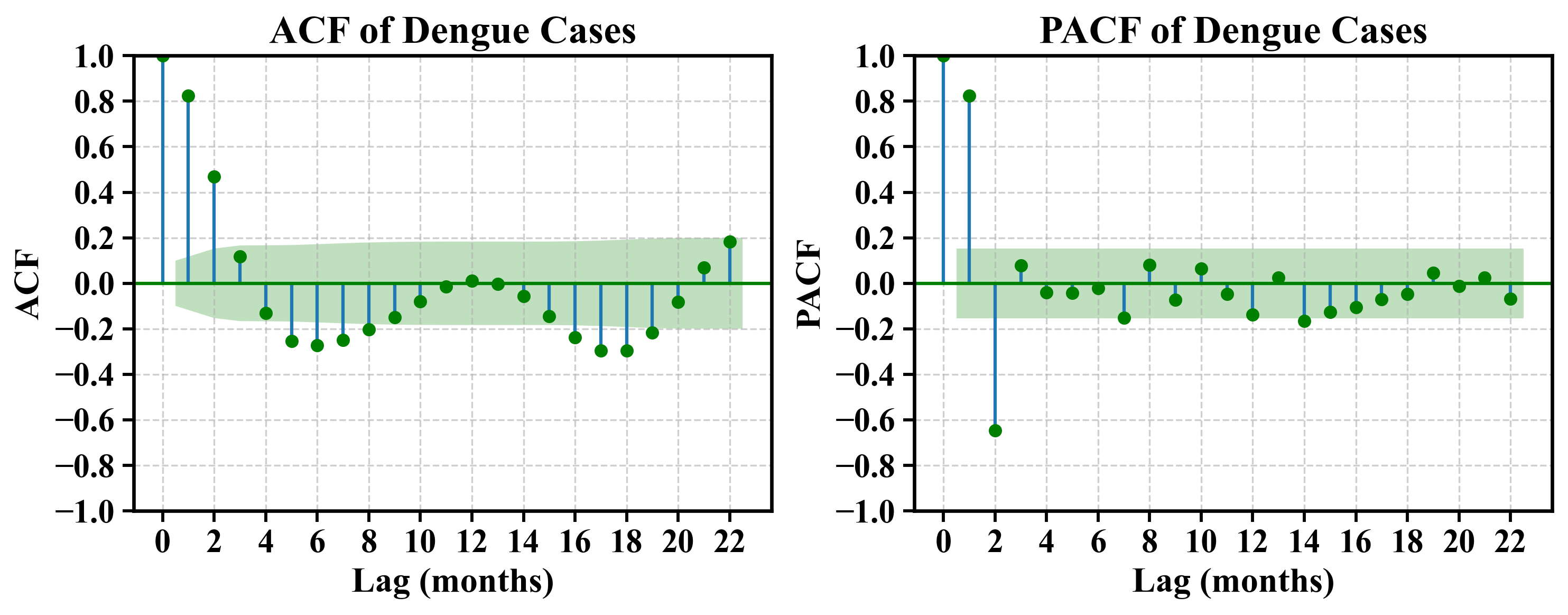}
	\caption{ACF and PACF plots of monthly dengue cases for Barishal.}	
	\label{acfb}
\end{figure}

\begin{figure}[H]
	\centering
	\subfloat[Dhaka]{\includegraphics[width=5.5 in]{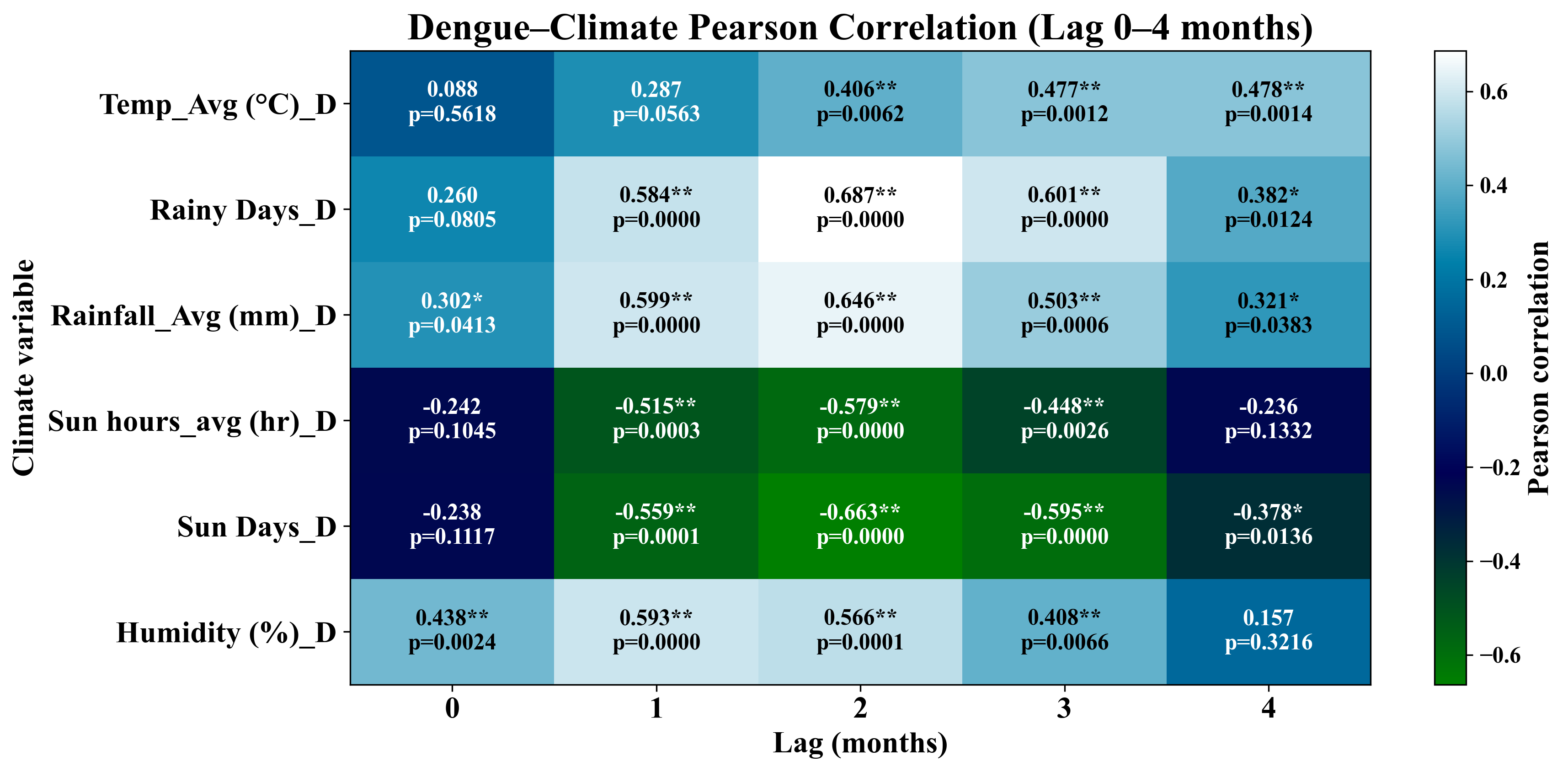}}\\
	\subfloat[Barishal]{\includegraphics[width=5.5 in]{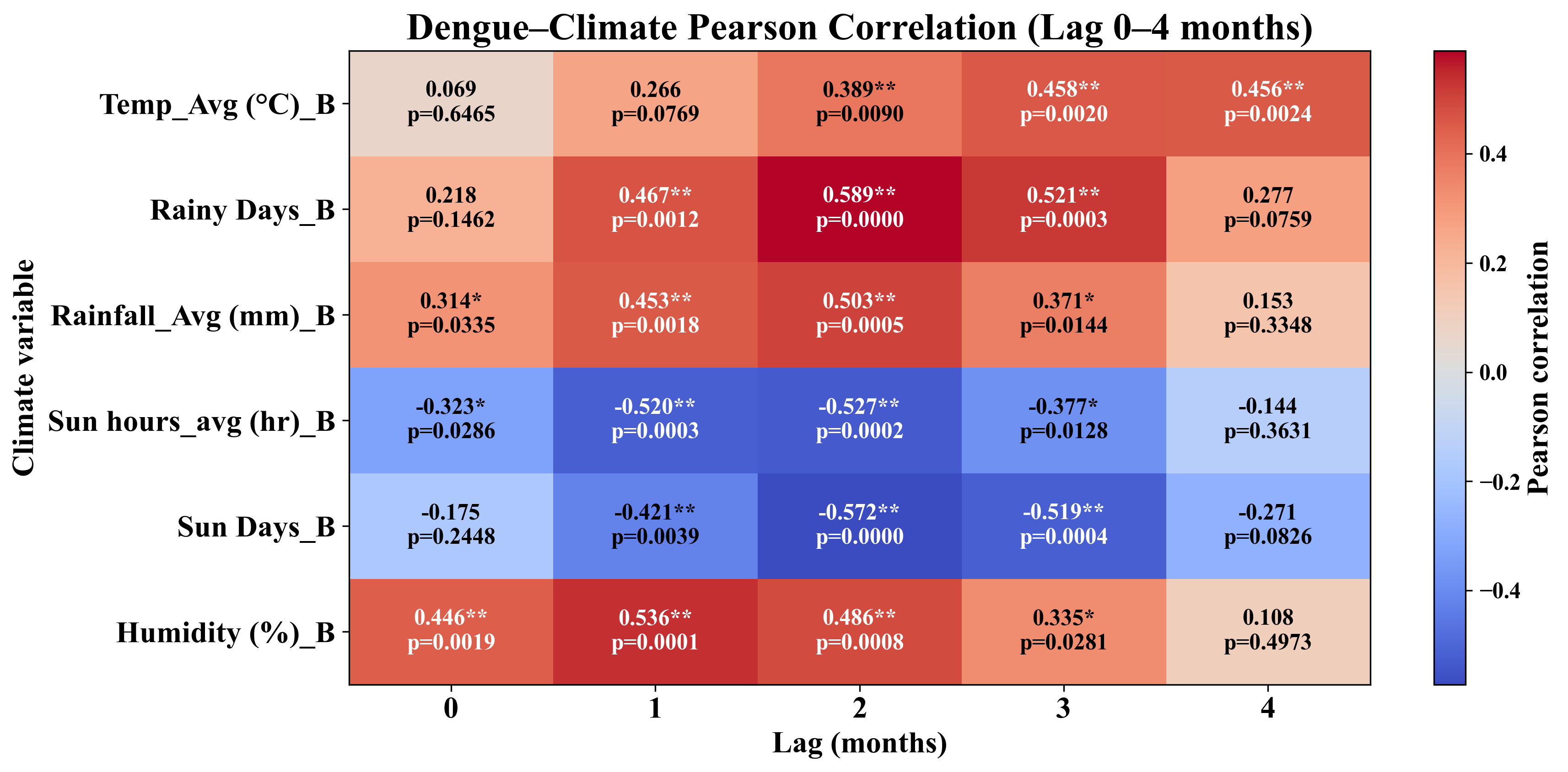}}
	\caption{Heatmap of lagged (0-4 months) Pearson correlations between dengue cases and climate variables.}	
	\label{corr}
\end{figure}

\begin{figure}[H]
	\centering
	\subfloat[Dhaka]{\includegraphics[width=4.5 in]{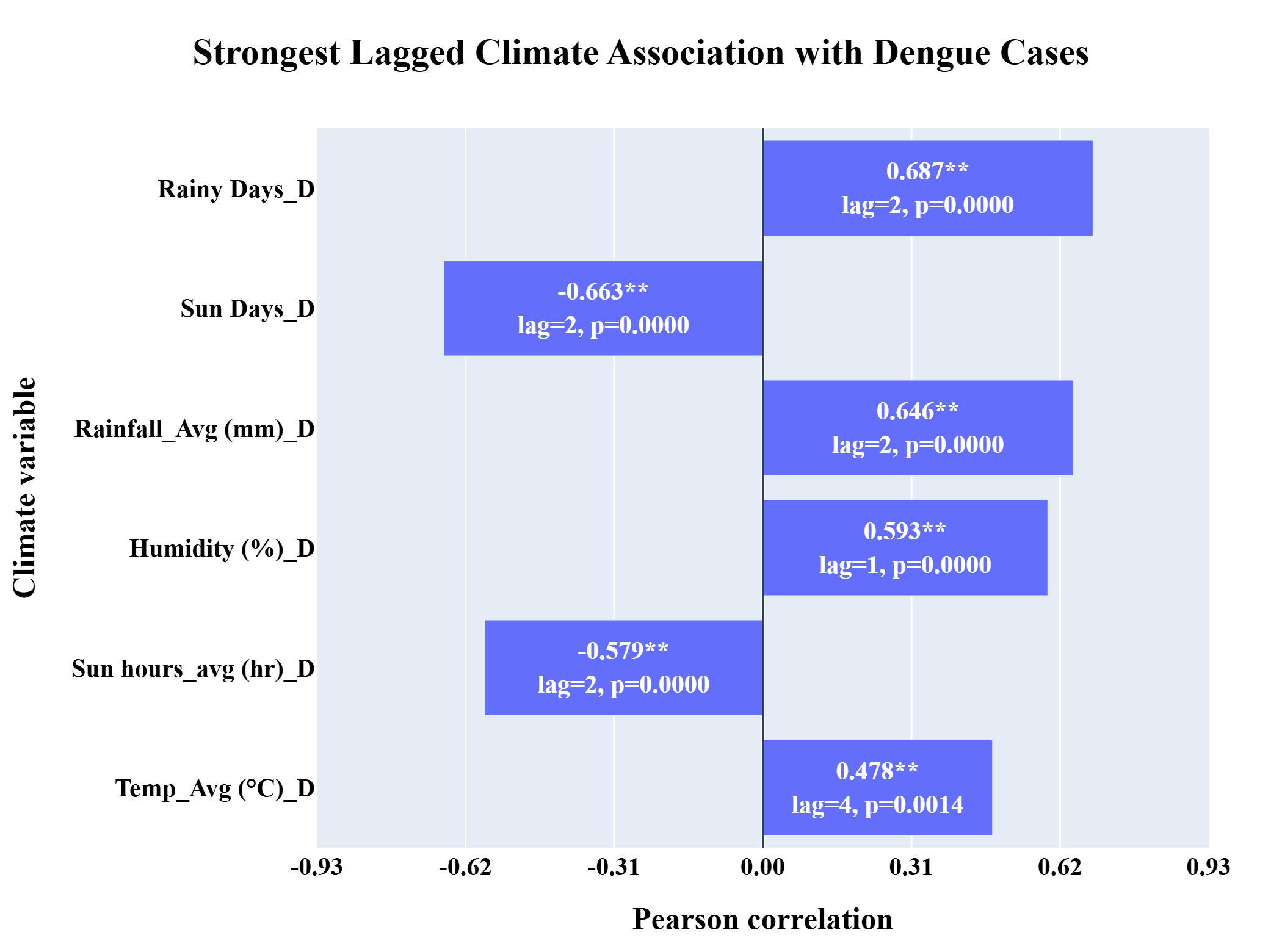}}\\
	\subfloat[Barishal]{\includegraphics[width=4.5 in]{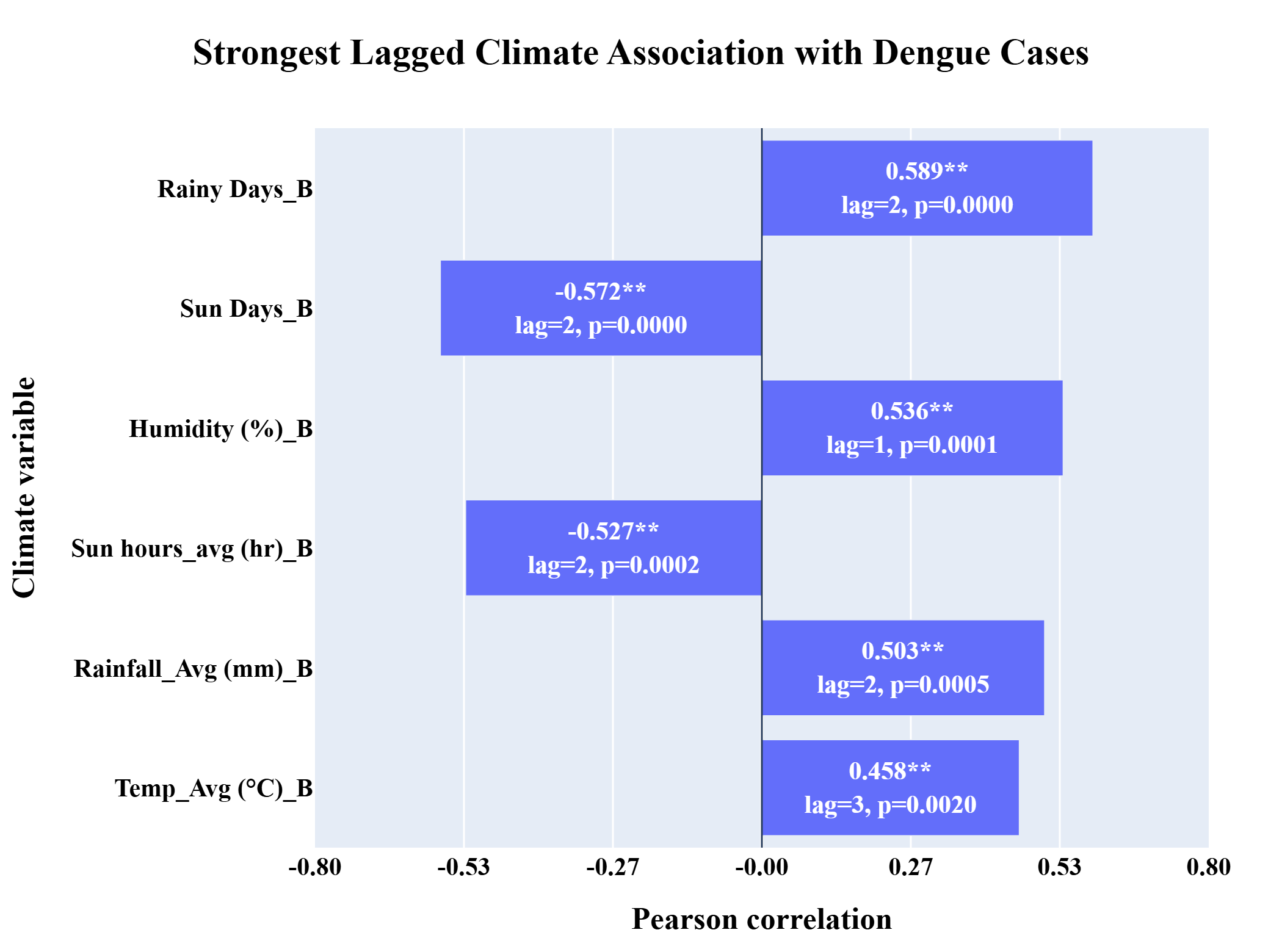}}
	\caption{Bar plot of the strongest lagged Pearson correlation for each climate variable.}	
	\label{bar}
\end{figure}

From thhe ACF-PACF analysis (Figure~\ref{acfd} and ~\ref{acfb}) show strong short-term persistence in both divisions: the ACF remains high at the first two lags ($\mathrm{ACF}(1)\approx 0.8$ and $\mathrm{ACF}(2)\approx 0.5$), and the PACF is dominated by lag 1 ($\mathrm{PACF}(1)\approx 0.8$) with an additional notable contribution at lag 2. These patterns support a low-order autoregressive specification, most naturally AR(1)-AR(2) (i.e. p=1 or 2), with a small moving-average component (e.g., q=1-2). In addition, the presence of a seasonal feature suggests an annual cycle in monthly dengue incidence, motivating the inclusion of seasonal terms with a 12-month period (s=12), such as a seasonal autoregressive component (P=1) and a seasonal moving average component (Q=1). 
\\
We compute lagged (0-4 months) Pearson correlations between dengue incidence and climate variables and summarize the strongest delayed effects using a heatmap (Figure~\ref{corr}(a-b)) and a ``best-lag" bar plot (Figure~\ref{bar}(a-b)). 

The results indicate that rainfall-related variables (rainy days and rainfall) show the strongest positive associations at about a 2-month lag in both Dhaka and Barishal, while humidity is most strongly positive at a 1-month lag. In contrast, sunshine-related measures (sun days and sunshine hours) are most strongly negative around a 2-month lag, and temperature shows a weaker but positive association at longer lags (3-4 months).
\subsection{SARIMAX}
To model monthly dengue incidence while accounting for annual seasonality and climate forcing, we apply Seasonal ARIMA with exogenous regressors (SARIMAX) to both Dhaka and Barishal using the same model-design framework and division wise dengue climate time series. Climate effects are incorporated through lagged predictors using a fixed lag structure selected from the strongest Pearson correlations: temperature (lag 3 months), rainy days (lag 2), rainfall (lag 2), sunshine hours (lag 2), sun days (lag 2), and humidity (lag 1).
We evaluate four exogenous feature sets: SET-1 (Temp + Rainy Days + Sun Hours + Humidity), SET-2 (Temp + Rainy Days + Sun Days + Humidity), SET-3 (Temp + Rainfall + Sun Days + Humidity), and SET-4 (Temp + Rainfall + Sun Hours + Humidity)--where the set definitions are the same across divisions but the exogenous values differ by division.
\begin{table}[H]
	
	\centering
	\setlength{\tabcolsep}{3.5pt}
	\caption{Best SARIMAX model per feature set-Dhaka}
	\label{sd}
	\resizebox{\textwidth}{!}{%}
	\begin{tabular}{l l c c c c c c}
		\hline
		\textbf{Feature set} & \textbf{Best SARIMAX} & \textbf{RMSE} & \textbf{MAE} & \textbf{MAPE  (\%)} & \textbf{AIC} & \textbf{BIC} & \textbf{Ljung--Box $p$} \\
		\hline
		SET-3 %(Temp+Rainfall+SunDays+Humidity)%
		 & SARIMAX(1,1,2)(1,0,0,12) & 2738.29 & 2444.24 & 231.84 & 457.60 & 467.42 & 0.45 \\
		SET-2 %(Temp+Rainy+SunDays+Humidity) 
		& SARIMAX(1,1,1)(1,0,0,12) & 5872.93 & 4974.29 & 372.01 & 449.81 & 458.54 & 0.02 \\
		SET-4 %(Temp+Rainfall+SunHours+Humidity) 
		& SARIMAX(1,1,2)(1,0,0,12) & 6972.58 & 5512.54 & 190.16 & 457.63 & 467.45 & 0.92 \\
		SET-1 %(Temp+Rainy+SunHours+Humidity) 
		& SARIMAX(2,1,2)(1,0,0,12) & 11300.41 & 9334.16 & 262.48 & 438.1 & 448.54 & 0.47 \\
		\hline
	\end{tabular}%
}
\end{table}

\begin{figure}[H]
	\centering 
	\subfloat[]{\includegraphics[width=5.5 in]{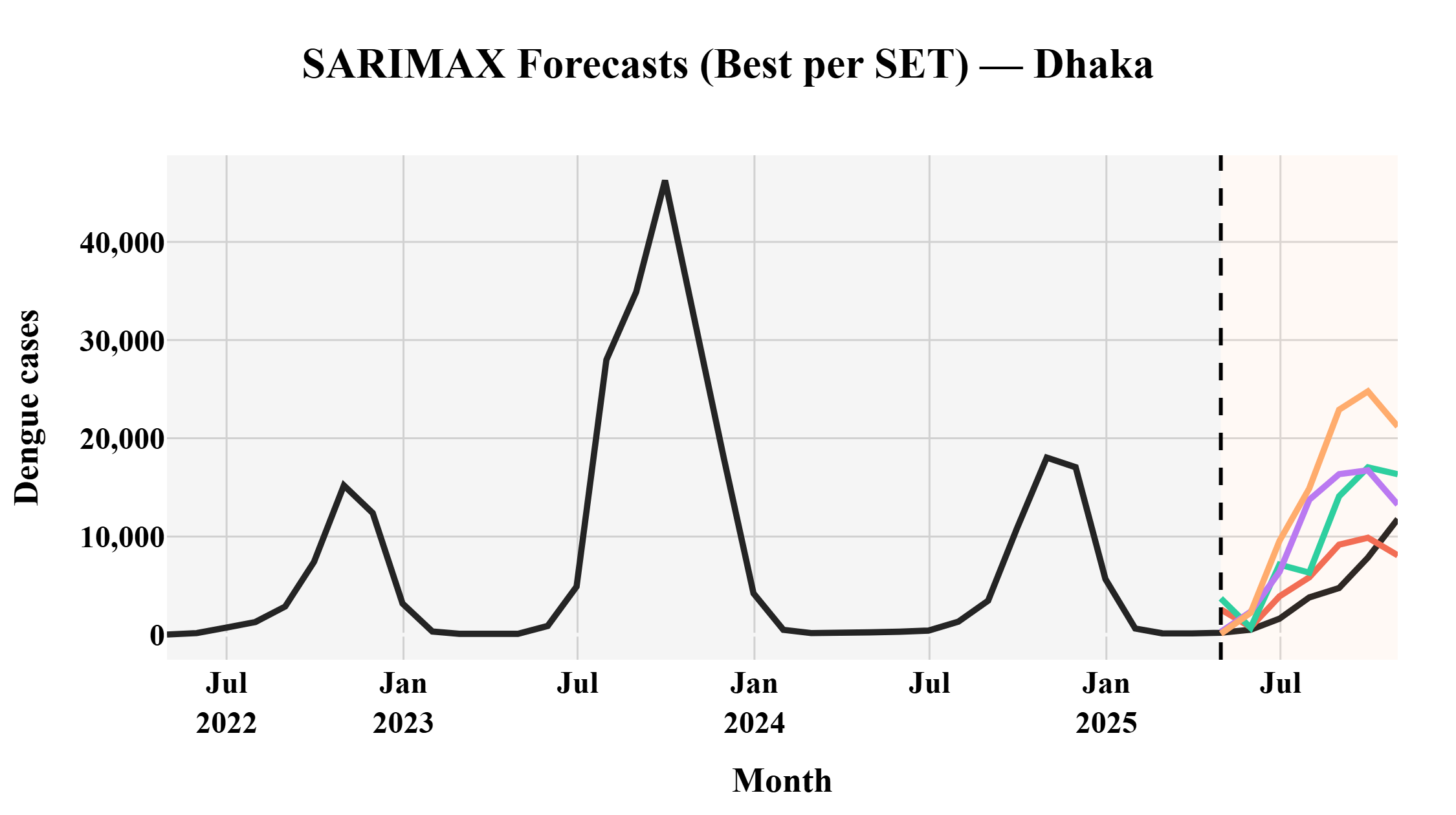}}\\
	\subfloat[]{\includegraphics[width=5.5 in]{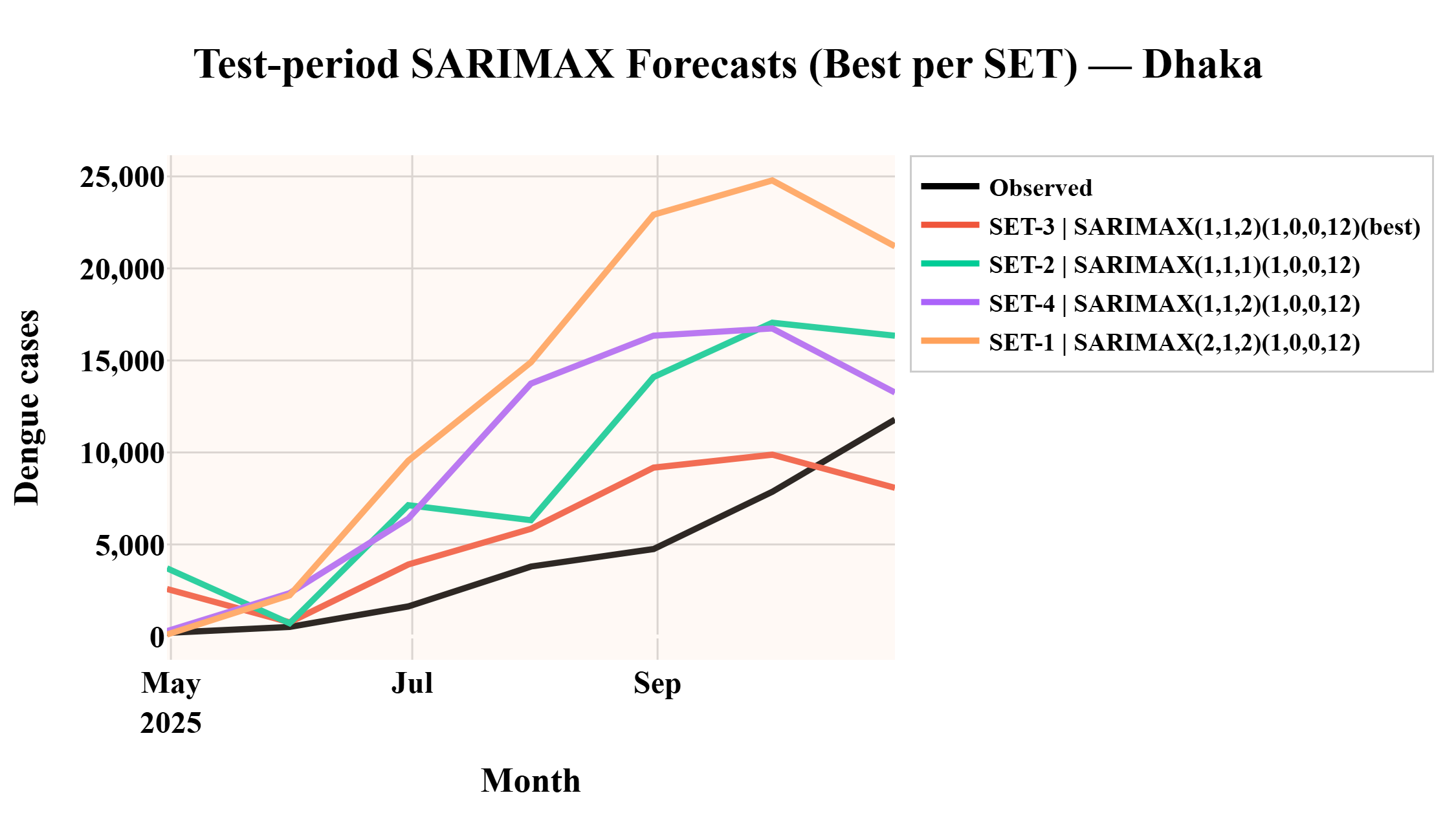}}
	\caption{SARIMAX results for Dhaka division. (a) Full series observed cases with test-window forecasts from the best model in each feature set; (b) Test period observed vs forecasts (best per feature set).}	
	\label{sard}
\end{figure}
\begin{figure}[H]
	\centering 
	\includegraphics[width=6 in]{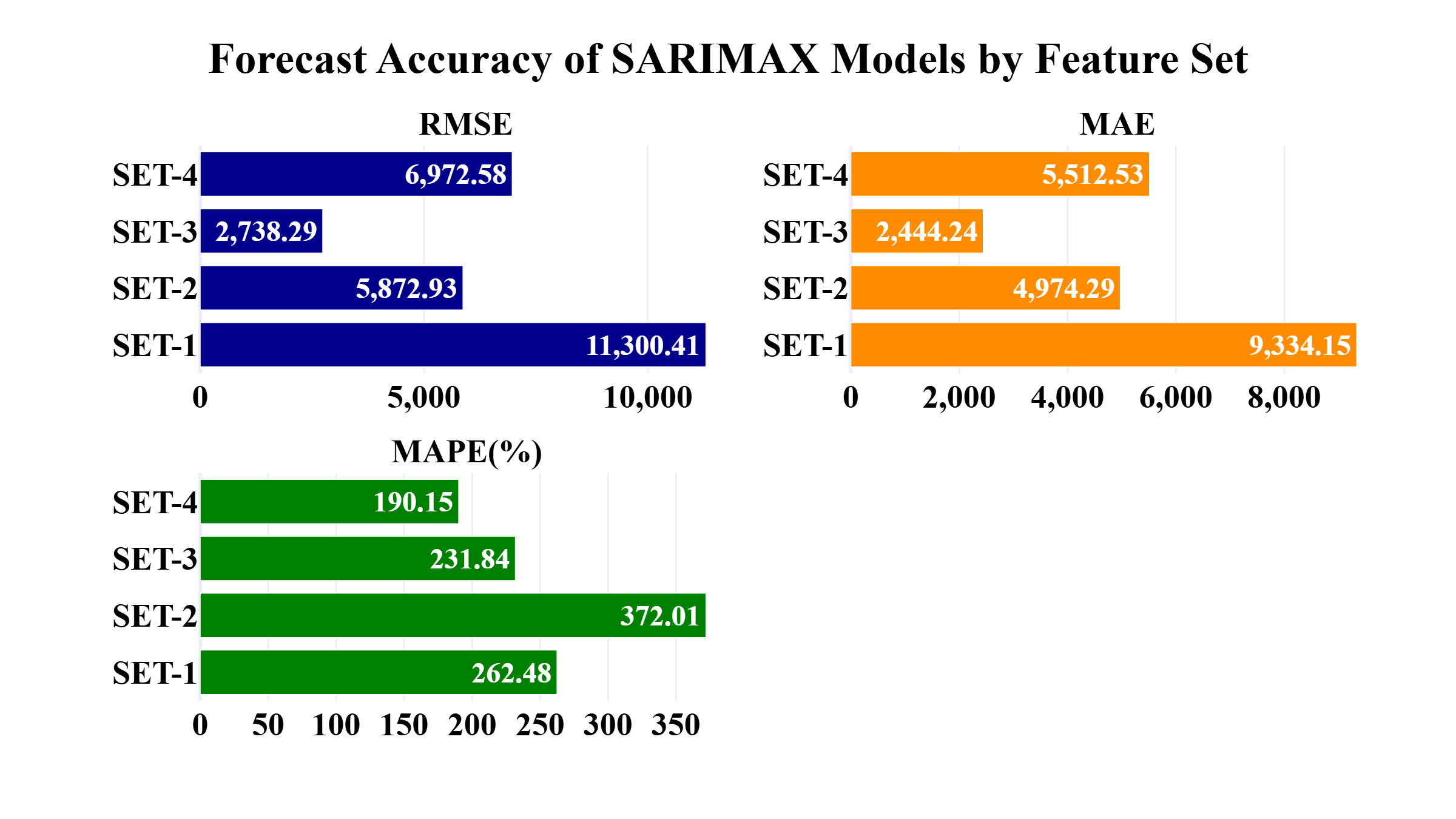}
	\caption{SARIMAX results for the Dhaka division, showing the error metrices (RMSE, MAE, MAPE) of the best performing model within each feature set.}	
	\label{sarde}
\end{figure}
For model structure, we evaluate a grid of seasonal SARIMAX candidates with annual seasonality (s=12), using orders $\text{SARIMAX}(p,d,q)\times(P, D, Q, 12)$  where $p, q\in\{0, 1, 2\}$, $d=1$, $P=1$, and $D, Q\in\{0, 1\}$. All candidates are fitted on the training data, forecasts are generated for the test period, and models are compared using forecast accuracy and residual diagnostics. For each feature set, we retain the specification with the lowest RMSE and report its AIC (Akaike Information Criterion), BIC (Bayesian Information Criterion), RMSE (Root Mean Squared Error), MAE (Mean Absolute Error), MAPE (\%)(Mean Absolute Percentage Error), and Ljung--Box $p$-value in Table~\ref{sd} and Table~\ref{sb}.

\begin{table}[H]
	
	\centering
	\setlength{\tabcolsep}{3.5pt}
	\caption{Best SARIMAX model per feature set-Barishal}
	\label{sb}
	\resizebox{\textwidth}{!}{%}
		\begin{tabular}{l l c c c c c c}
			\hline
			\textbf{Feature set} & \textbf{Best SARIMAX} & \textbf{RMSE} & \textbf{MAE} & \textbf{MAPE  (\%)} & \textbf{AIC} & \textbf{BIC} & \textbf{Ljung--Box $p$} \\
			\hline
			SET-2 %(Temp+Rainy+SunDays+Humidity)
			 & SARIMAX(0,1,1)(1,0,0,12) & 817.56 & 717.78 & 39.96 & 401.50 & 409.45 & 0.32 \\
			SET-1 %(Temp+Rainy+SunHours+Humidity)
			 & SARIMAX(0,1,1)(1,0,0,12) & 1186.95 & 888.33 & 56.29 & 399.97 & 407.92 & 0.50 \\
			SET-4 %(Temp+Rainfall+SunHours+Humidity)
			 & SARIMAX(1,1,0)(1,0,0,12) & 1306.66 & 967.98 & 53.72 & 391.42 & 399.06 & 0.44 \\
			SET-3 %(Temp+Rainfall+SunDays+Humidity) 
			& SARIMAX(2,1,1)(1,1,1,12) & 2026.49 & 1379.17 & 56.76 & 165.59 & 167.56 & 0.01 \\
			\hline
		\end{tabular}%
	}
\end{table}
Tables~\ref{sd} and \ref{sb} summarize the best SARIMAX specification for each feature set in Dhaka and Barishal based on forecast accuracy and residual diagnostics. The Ljung--Box $p$-value is used to assess whether residual autocorrelation remains, with higher $p$-values indicating better model adequacy. The results show that SET-3 gives the most favorable balance of accuracy and residual behavior in Dhaka, while SET-2 does so in Barishal, highlighting division-specific differences in the usefulness of the climate feature sets. 
\begin{figure}[H]
	\centering 
	\subfloat[]{\includegraphics[width=5.5 in]{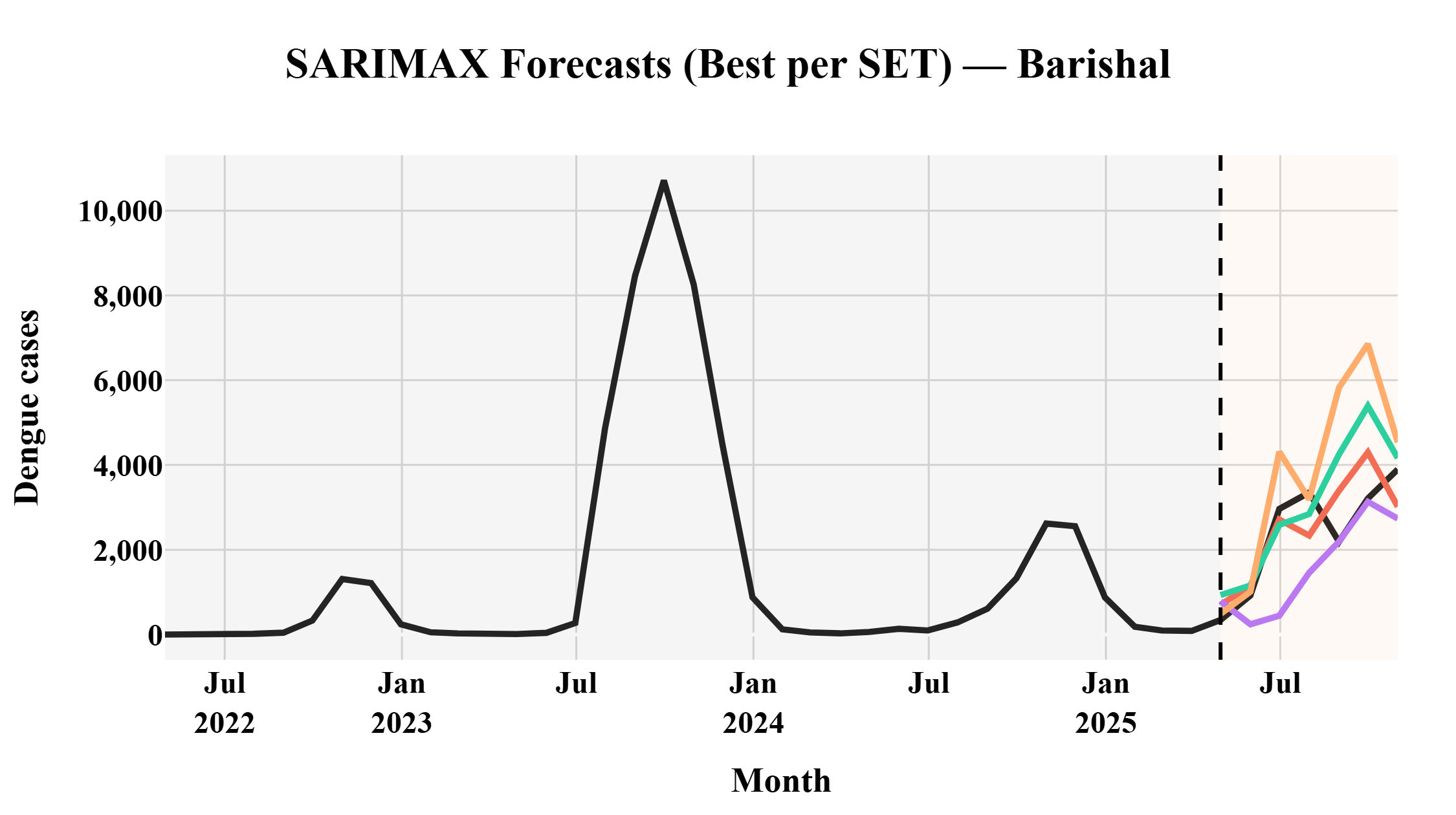}}\\
	\subfloat[]{\includegraphics[width=5.5 in]{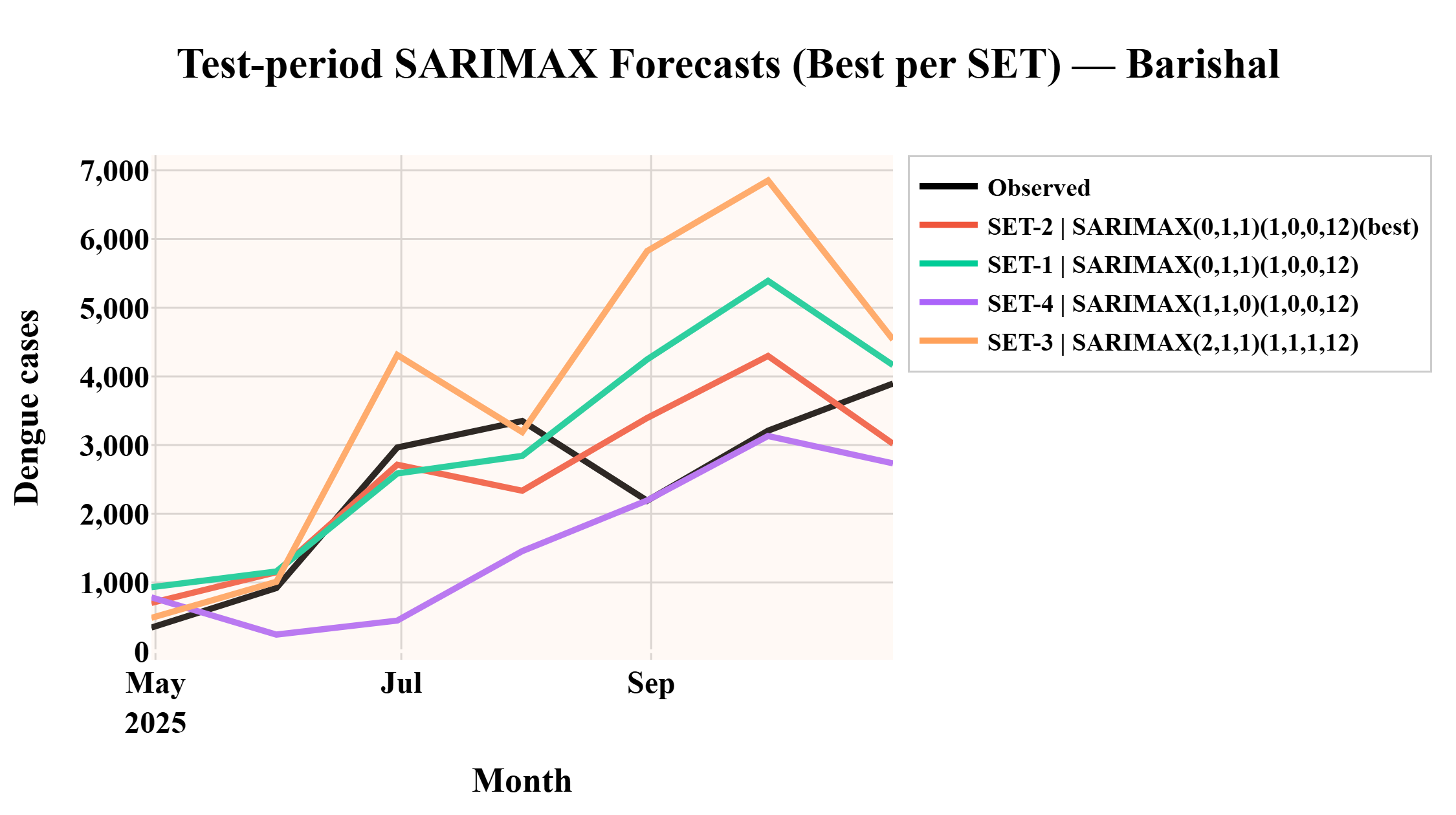}}
	\caption{SARIMAX results for Barishal division. (a) Full series observed cases with test-window forecasts from the best model in each feature set; (b) Test period observed vs forecasts (best per feature set).}	
	\label{sarb}
\end{figure}
\begin{figure}[H]
	\centering 
	\includegraphics[width=6 in]{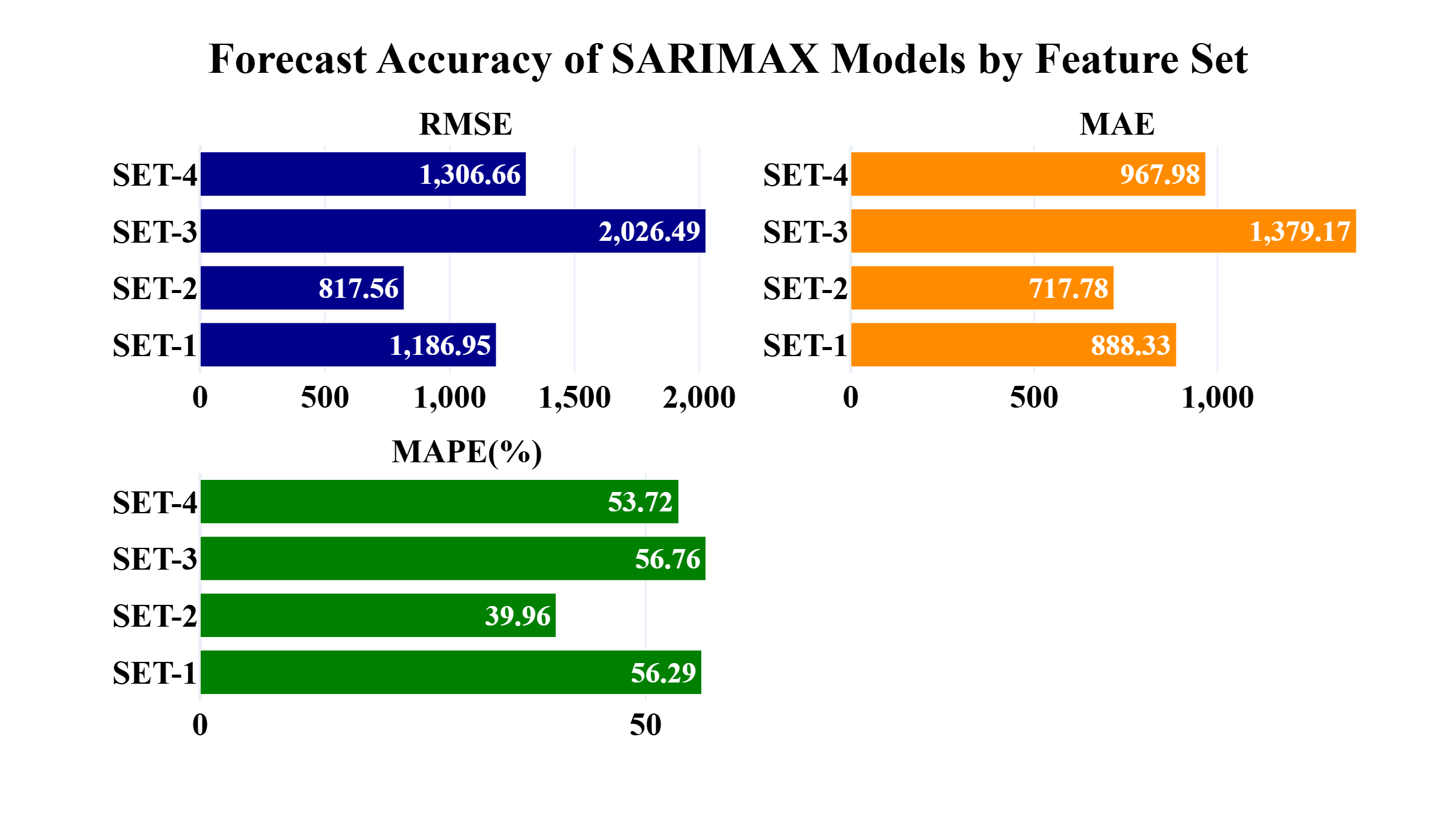}
	\caption{SARIMAX results for the Barishal division, showing the error metrices (RMSE, MAE, MAPE) of the best performing model within each feature set.}	
	\label{sarbe}
\end{figure}
Based on this selection, the overall best model is SET-3 with SARIMAX(1, 1, 2)(1, 0, 0, 12) for Dhaka (Figure~\ref{sard}(a-b)) (RMSE=2738.29) (Figure~\ref{sarde}) and  SET-2 with SARIMAX(0, 1, 1)(1, 0, 0, 12) for Barishal (Figure~\ref{sarb}(a-b)) (RMSE= 817.56)(Figure~\ref{sarbe}); the train-test forecast plots and the metric bar charts summarize these set-wise comparisons for each division. 
\subsection{MPR}
\begin{figure}[H]
	\centering 
	\subfloat[]{\includegraphics[width=4.2 in]{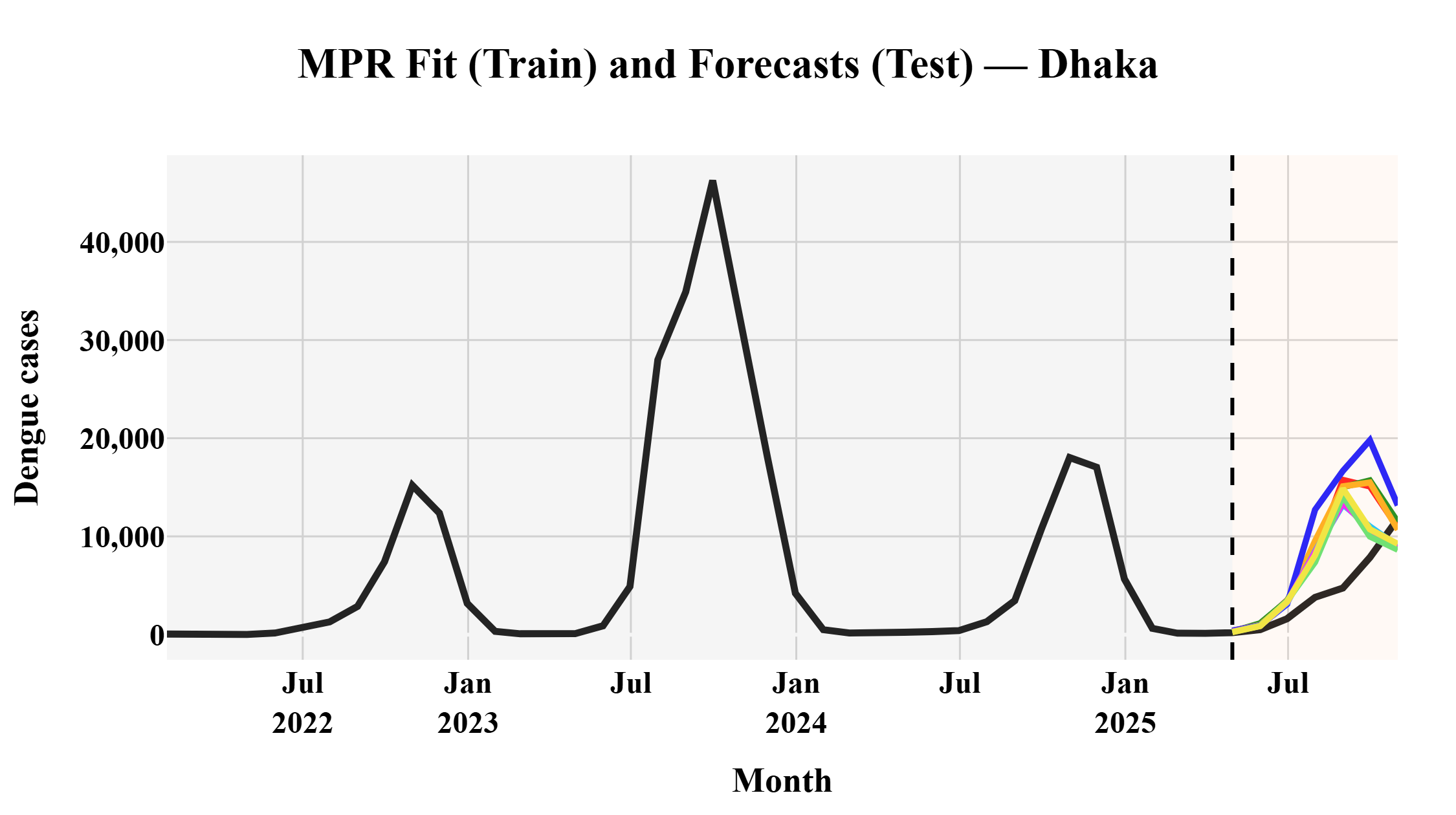}}\\
	\subfloat[]{\includegraphics[width=4.2 in]{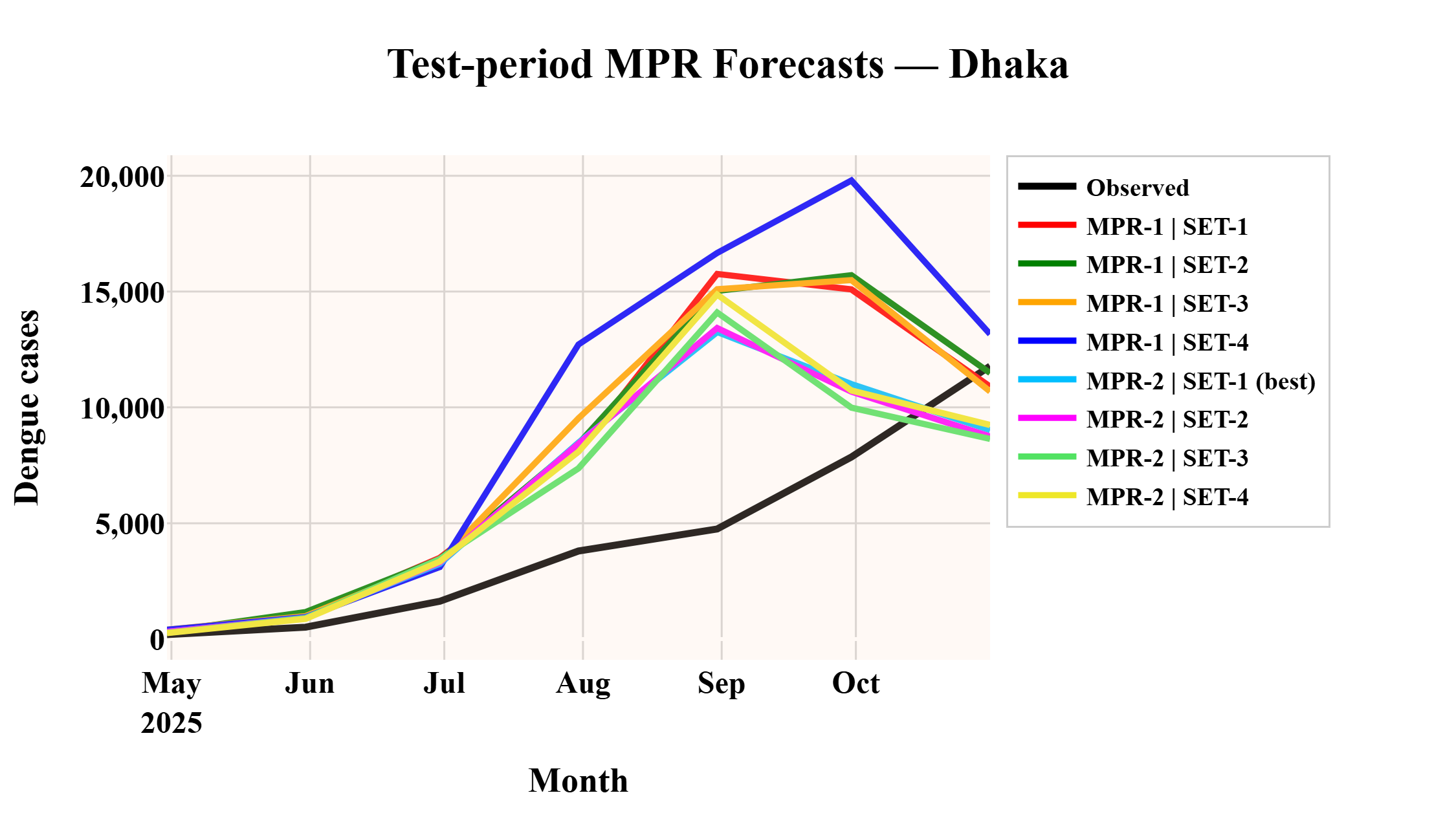}}
	
	\caption{MPR results for Dhaka division. (a) Full series observed cases with test-window forecasts from the best model in each feature set; (b) Test period observed vs forecasts (best per feature set).}	
	\label{mprd}
\end{figure}

\begin{figure}[H]
	\centering 
	
	\includegraphics[width=6 in]{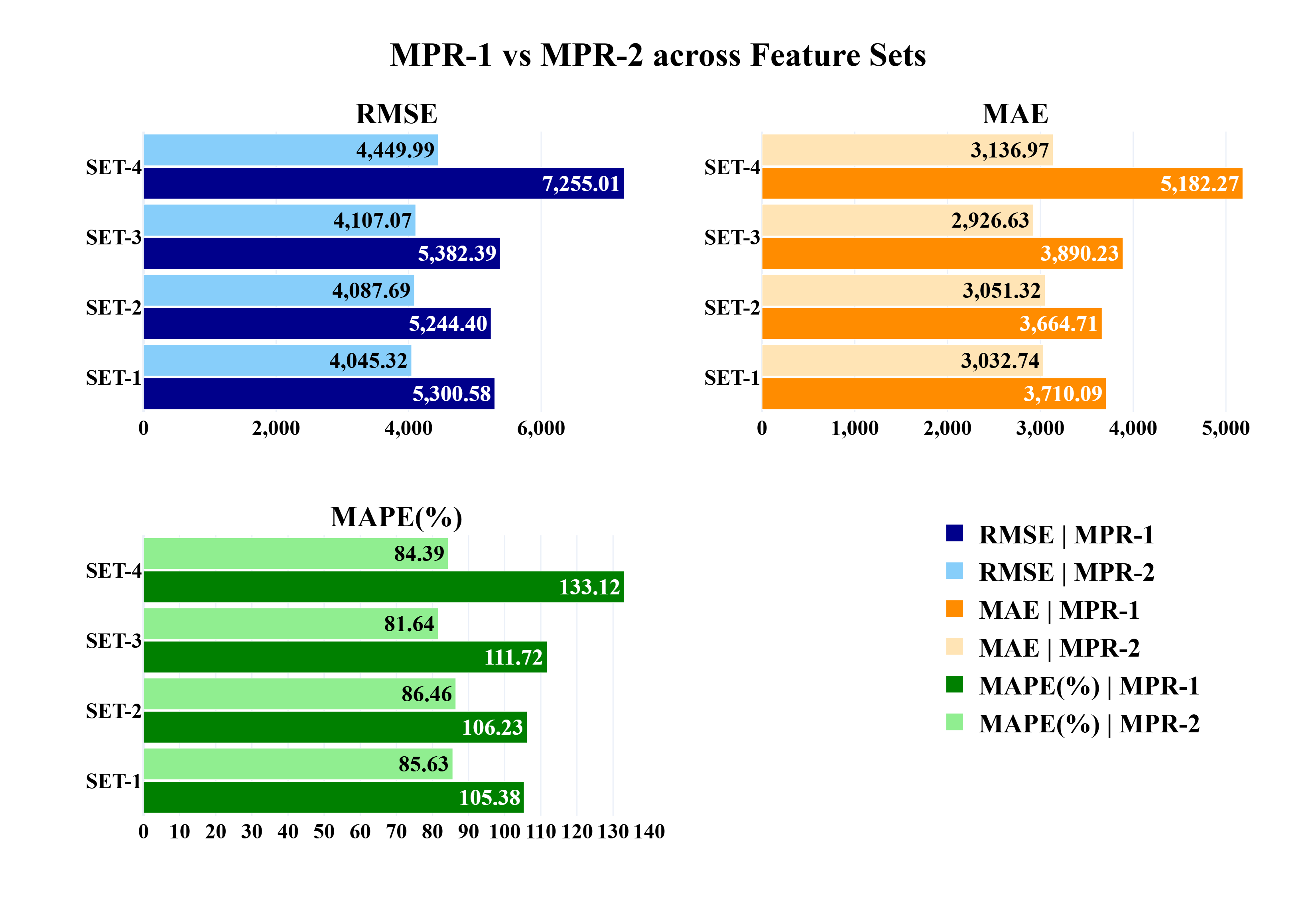}
	\caption{MPR results for the Dhaka division, showing the error metrices (RMSE, MAE, MAPE) of the best performing model within each feature set.}	
	\label{mprde}
\end{figure}

\begin{figure}[H]
	\centering 
	\subfloat[]{\includegraphics[width=4.2 in]{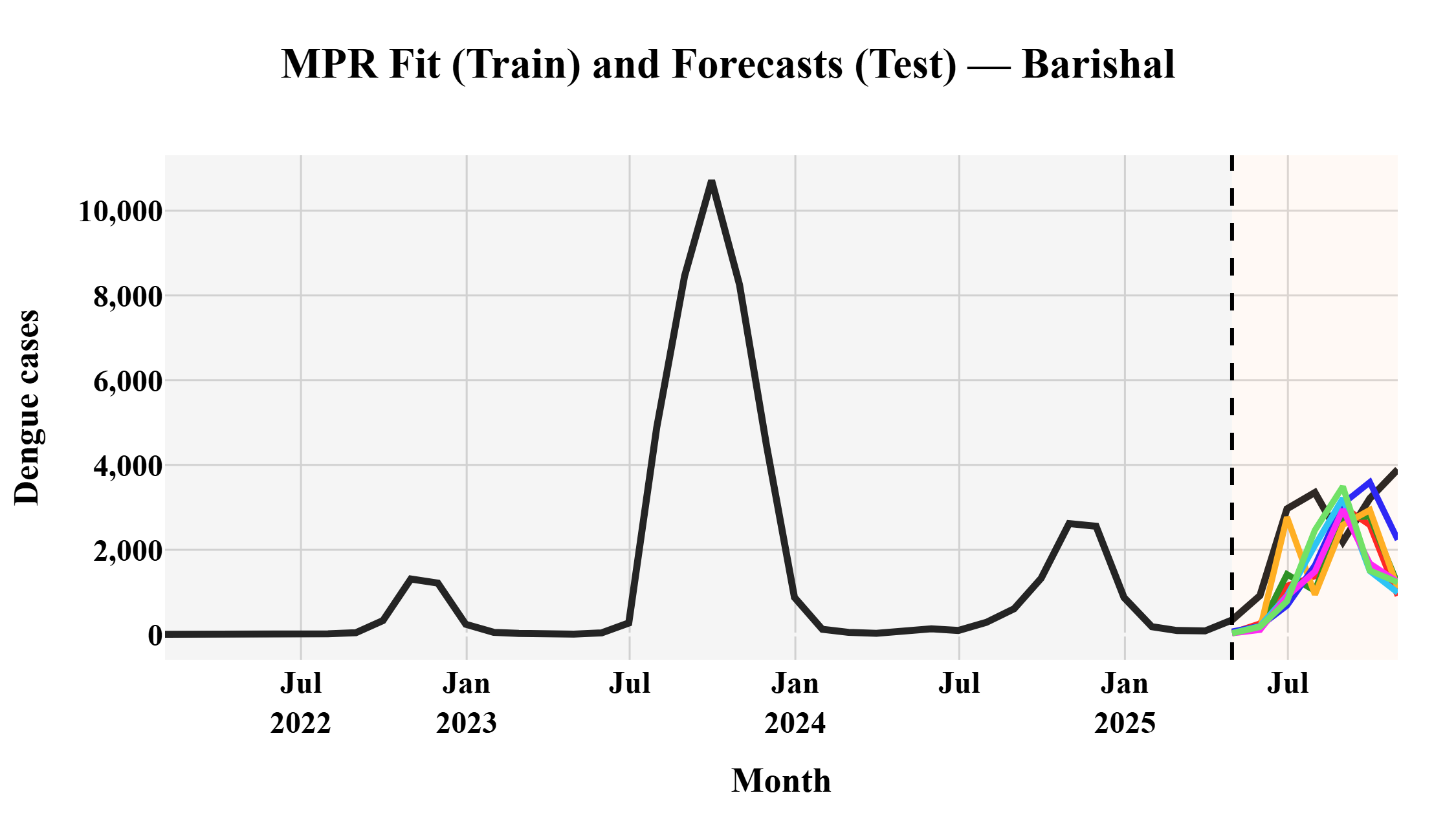}}\\
	\subfloat[]{\includegraphics[width=4.2 in]{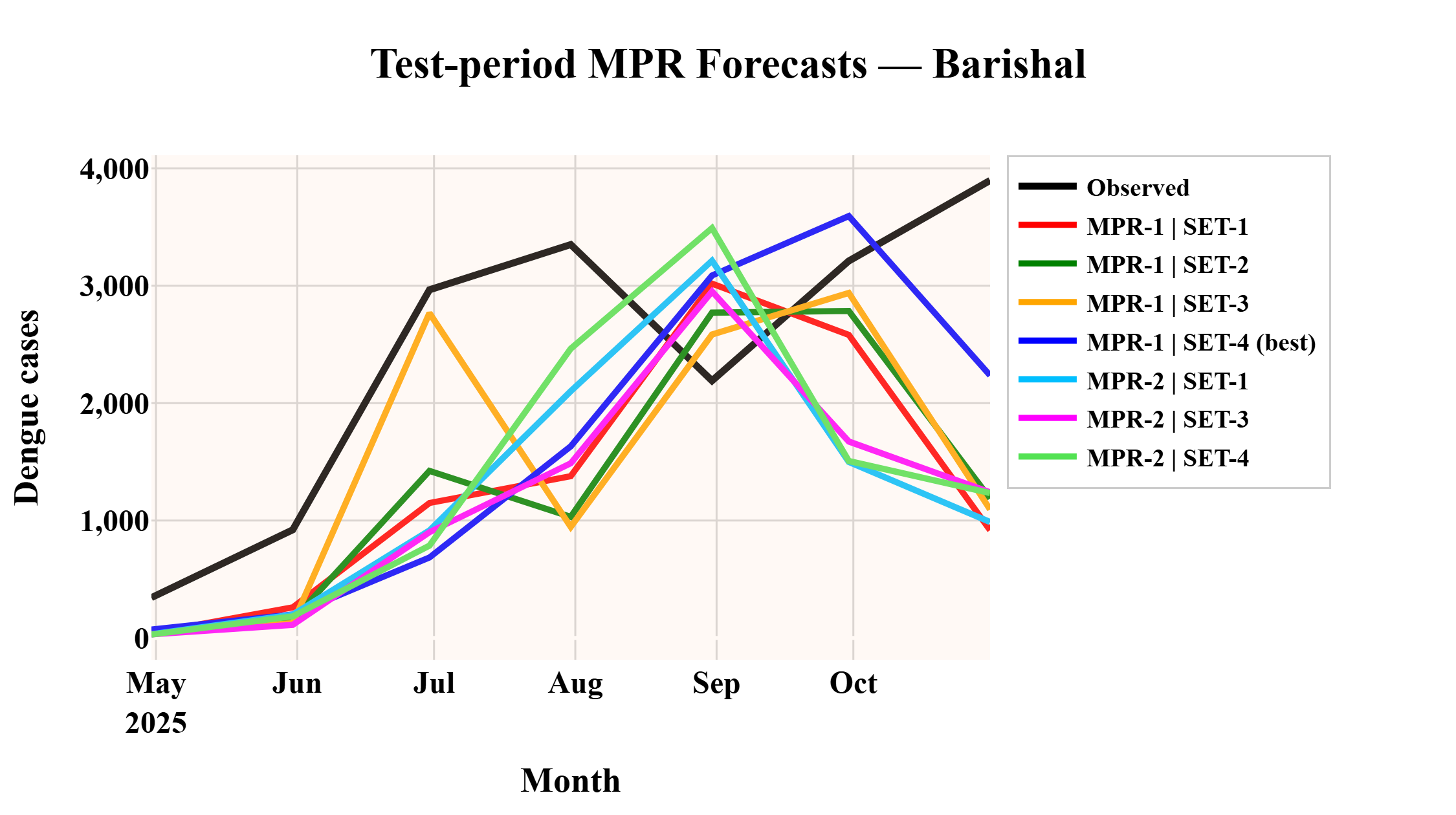}}
	
	\caption{MPR results for Barishal division. (a) Full series observed cases with test-window forecasts from the best model in each feature set; (b) Test period observed vs forecasts (best per feature set).}	
	\label{mprb}
\end{figure}

\begin{figure}[H]
	\centering 
	
	\includegraphics[width=6 in]{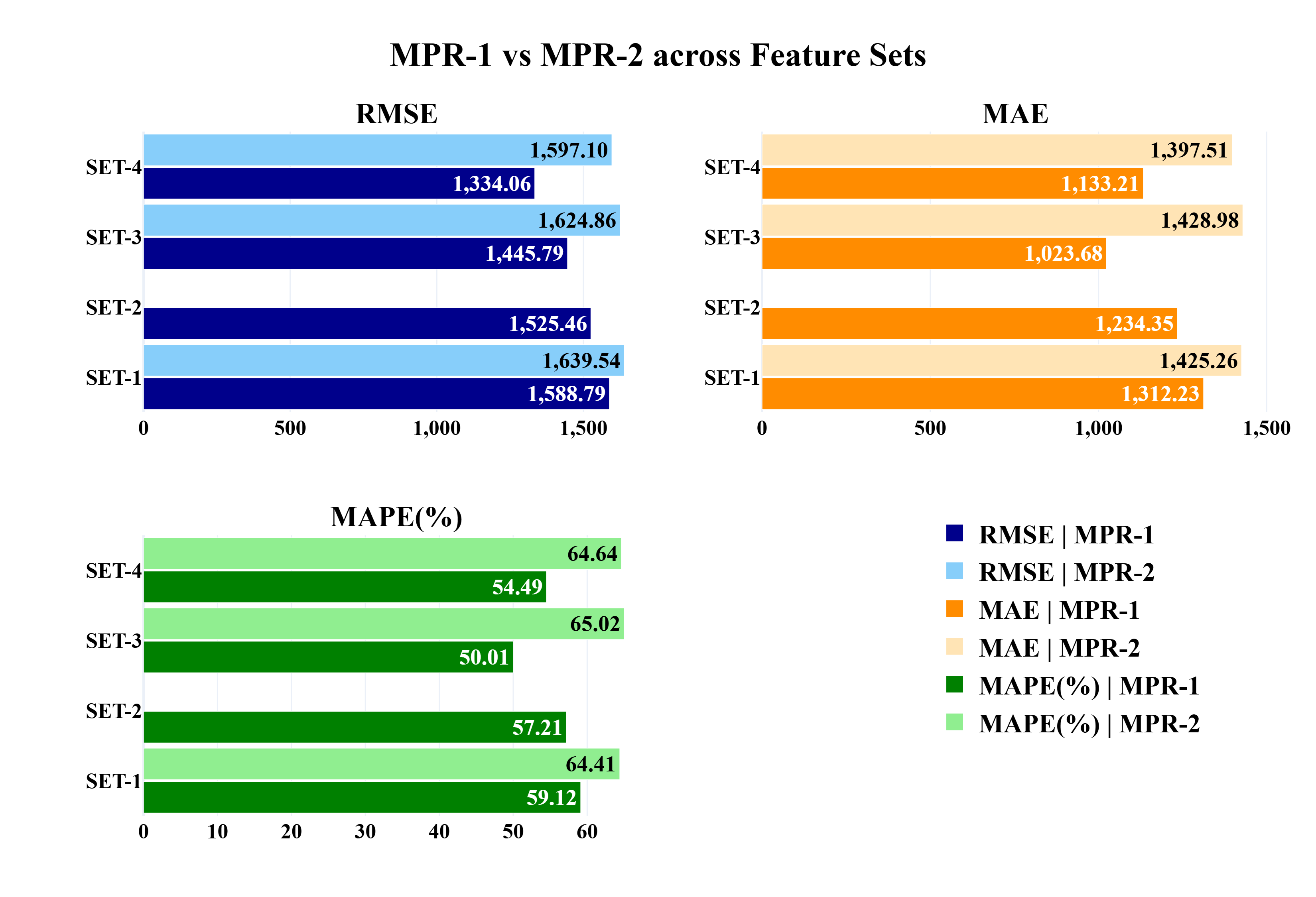}
	\caption{MPR results for the Barishal division, showing the error metrices (RMSE, MAE, MAPE) of the best performing model within each feature set.}	
	\label{mprbe}
\end{figure}

Across the MPR experiments, the model ranking differs by division. In Dhaka (Figure~\ref{mprde}), the lowest test error is obtained by MPR-2 with SET-1 (temperature, rainy days, sunshine hours, humidity) with RMSE=4045.32, MAE=3032.74, and MAPE=85.63$\%$. This indicates that adding the lagged dengue term (lag 1) (short-term persistence) improves test-period prediction relative to the climate-only MPR-1 runs. In this best Dhaka specification (Figure~\ref{mprd}(a-b)), dengue variability aligns most consistently with temperature and humidity (mosquito development and survival), rainy days (repeated refilling of breeding sites), and sunshine hours (drying and habitat persistence), while the case-lag captures the momentum from the most recent outbreak level.
In Barishal (Figure~\ref{mprbe}), the best-performing configuration is SET-4(temperature, rainfall, sunshine hours, humidity), where MPR-1 is overall best with RMSE=1334.06, MAE=1133.21, and MAPE=54.49$\%$; the best case-lag (lag 1) variant is MPR-2 with SET-4 (RMSE=1597.10, MAE=1397.51, and MAPE=64.64$\%$). This suggests Barishal's dengue dynamics are comparatively more explainable by climate forcing alone (Figure~\ref{mprb}(a-b)), especially rainfall amount (water availability for larval habitats), together with temperature and humidity, while sunshine hours likely reflects drying and standing water persistence that modulates how effectively rainfall translates into breeding suitability. 
Additionally, MPR-2 with SET-2 in Barishal produces extremely large test errors (RMSE=159493.23, MAE=61618.84, MAPE=2095.06$\%$), indicating an unrealistic forecast for this specification; therefore, it is excluded from the time-series forecast figure (Figure~\ref{mprb}(a-b)) and the performance bar plot (Figure~\ref{mprbe}) for Barishal to avoid compressing the scale and obscuring differences among the remaining models.
\subsection{ANN}
\begin{figure}[H]
	\centering 
	\subfloat[]{\includegraphics[width=4.2 in]{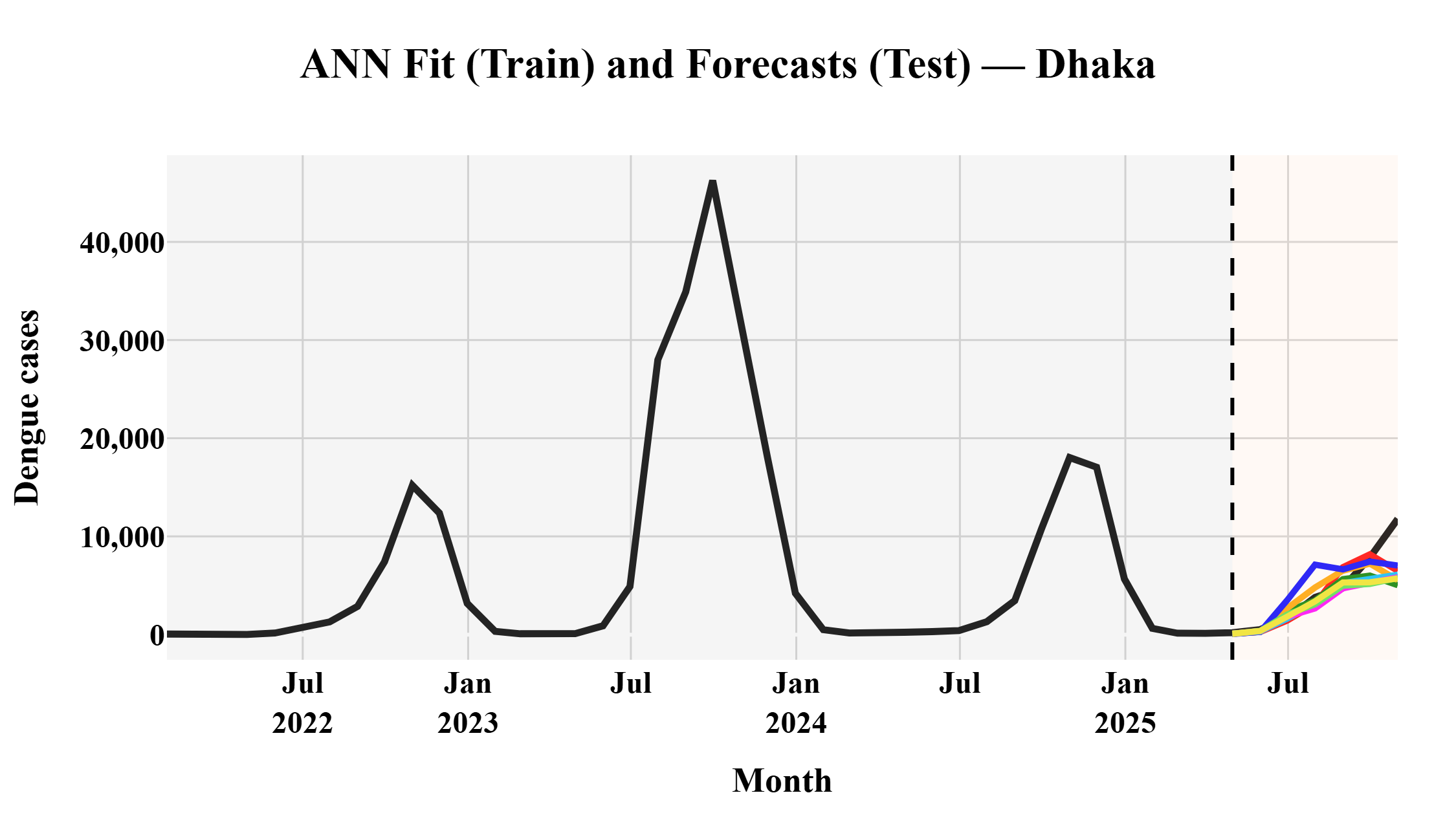}}\\
	\subfloat[]{\includegraphics[width=4.2 in]{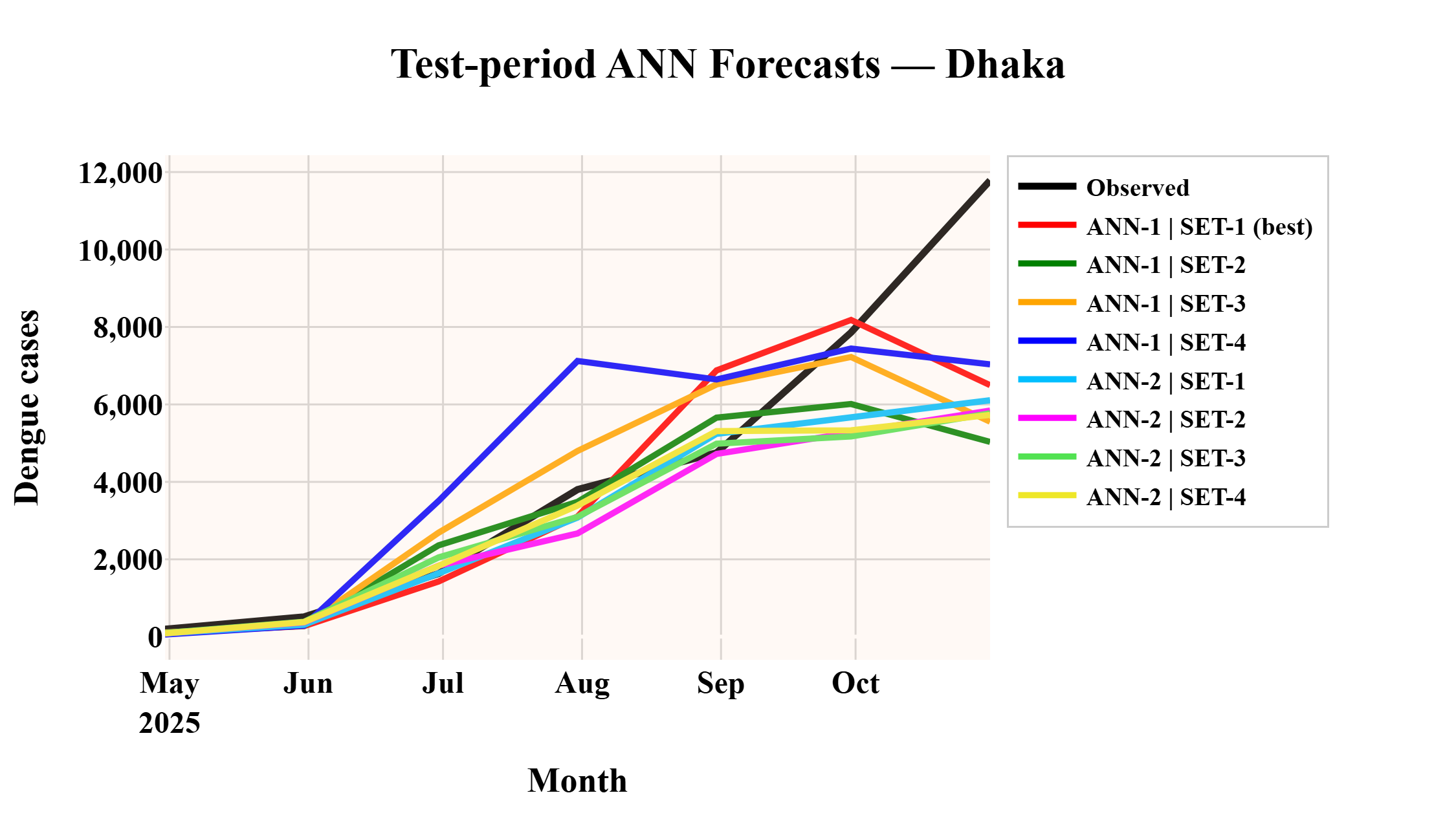}}
	
	\caption{ANN results for Dhaka division. (a) Full series observed cases with test-window forecasts from the best model in each feature set; (b) Test period observed vs forecasts (best per feature set).}	
	\label{annd}
\end{figure}

\begin{figure}[H]
	\centering 
	
	\includegraphics[width=6 in]{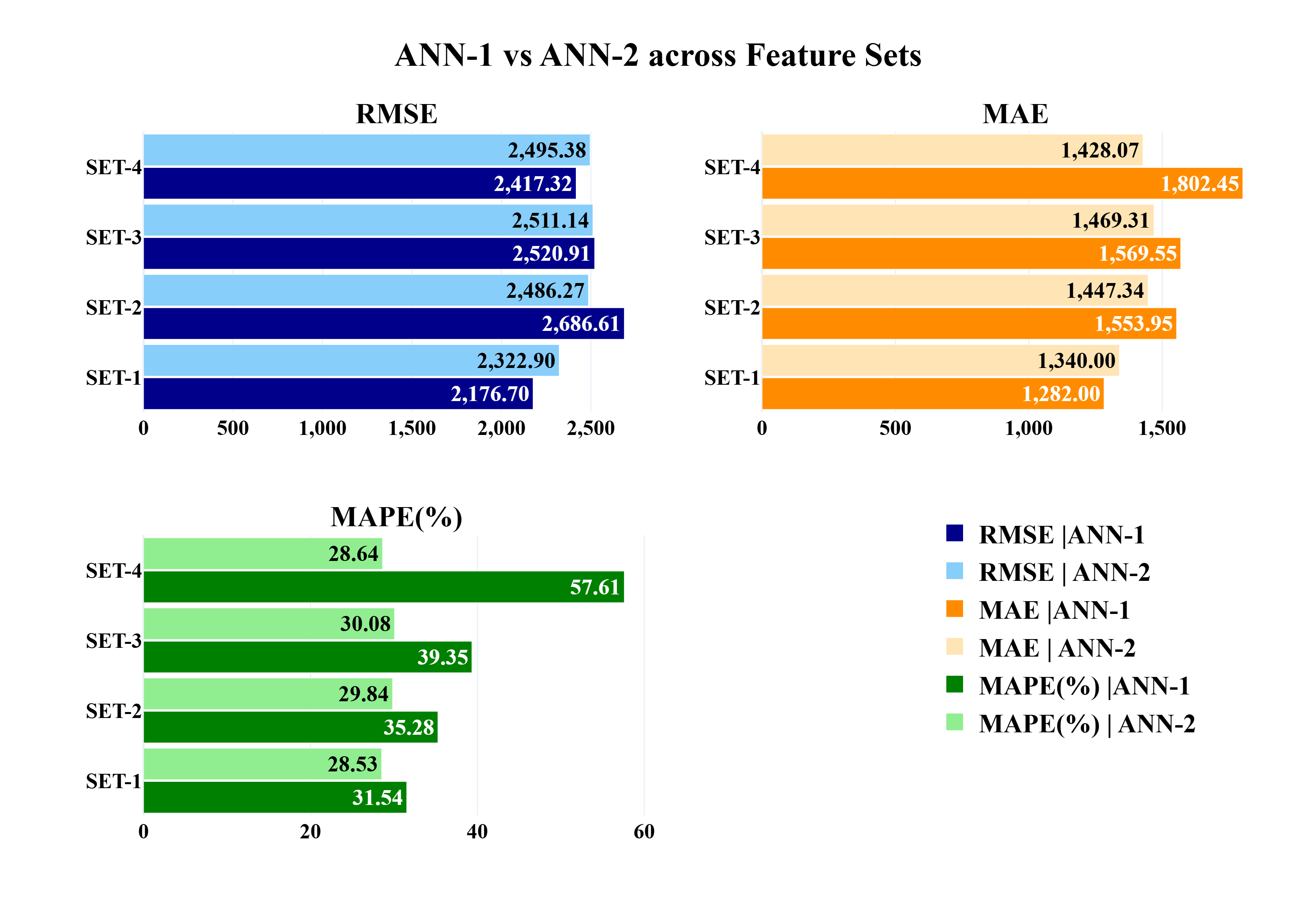}
	\caption{ANN results for the Dhaka division, showing the error metrices (RMSE, MAE, MAPE) of the best performing model within each feature set.}	
	\label{annde}
\end{figure}

\begin{figure}[H]
	\centering 
	\subfloat[]{\includegraphics[width=4.2 in]{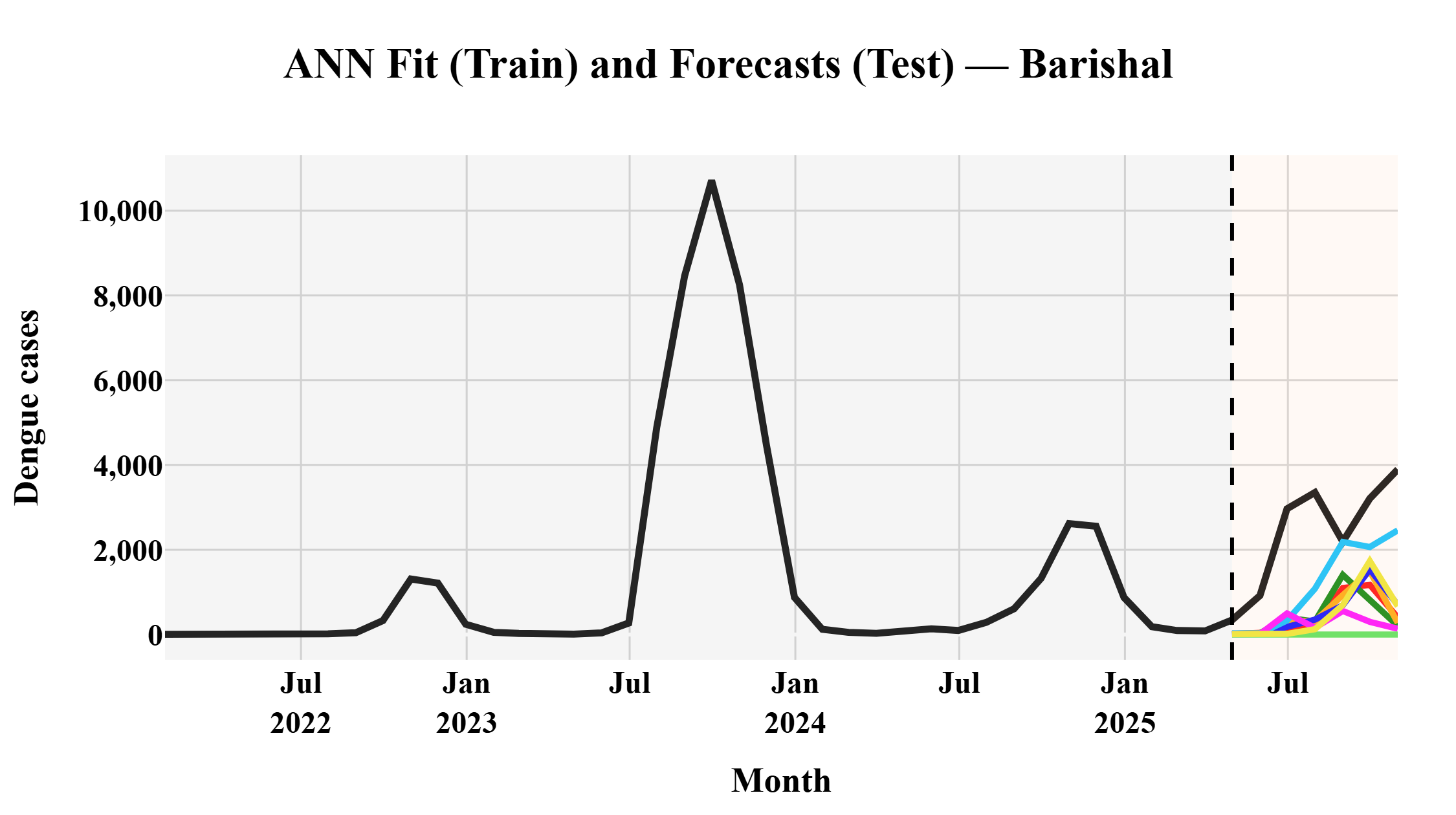}}\\
	\subfloat[]{\includegraphics[width=4.2 in]{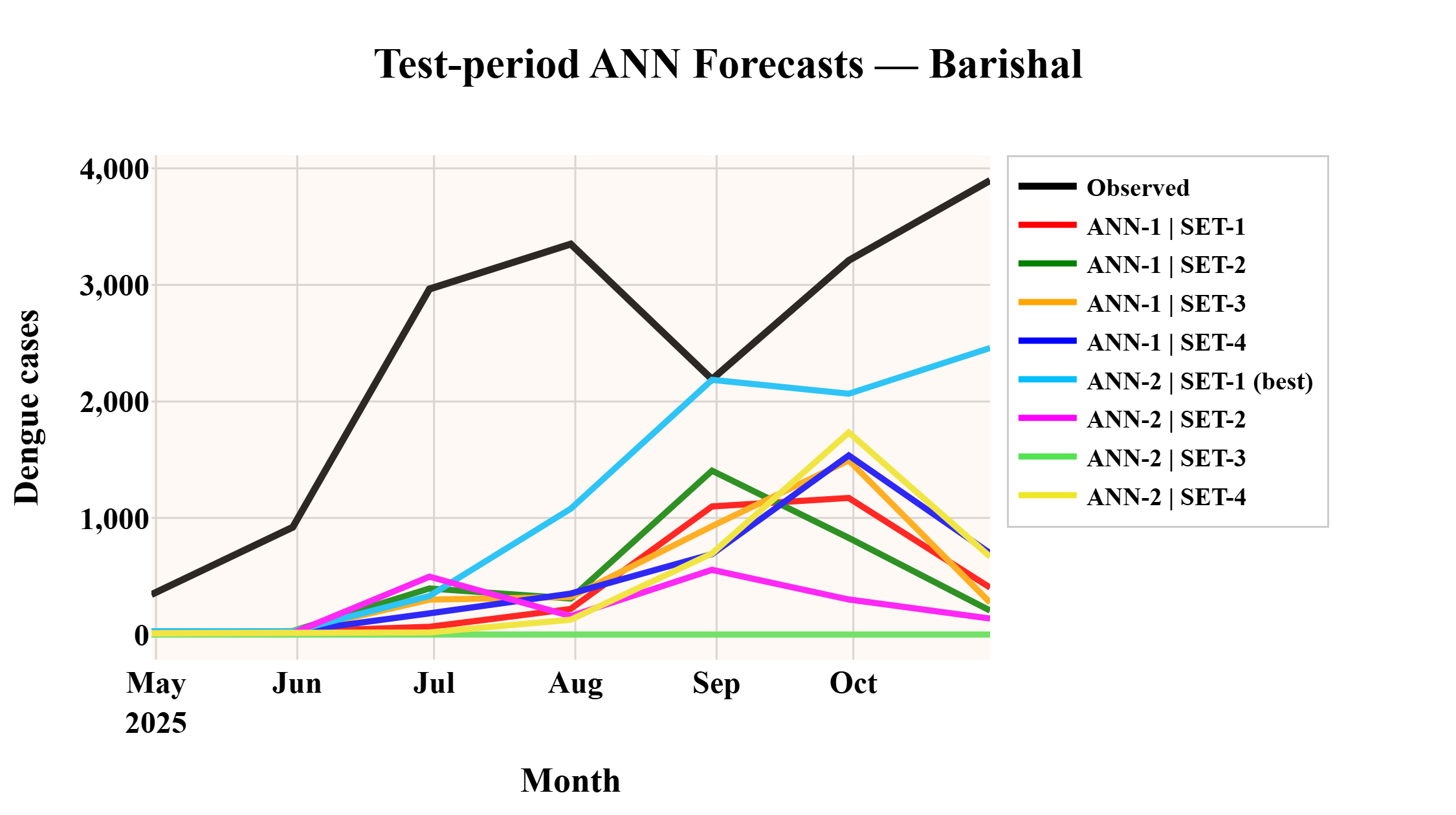}}
	
	\caption{ANN results for Barishal division. (a) Full series observed cases with test-window forecasts from the best model in each feature set; (b) Test period observed vs forecasts (best per feature set).}	
	\label{annb}
\end{figure}

\begin{figure}[H]
	\centering 
	
	\includegraphics[width=6 in]{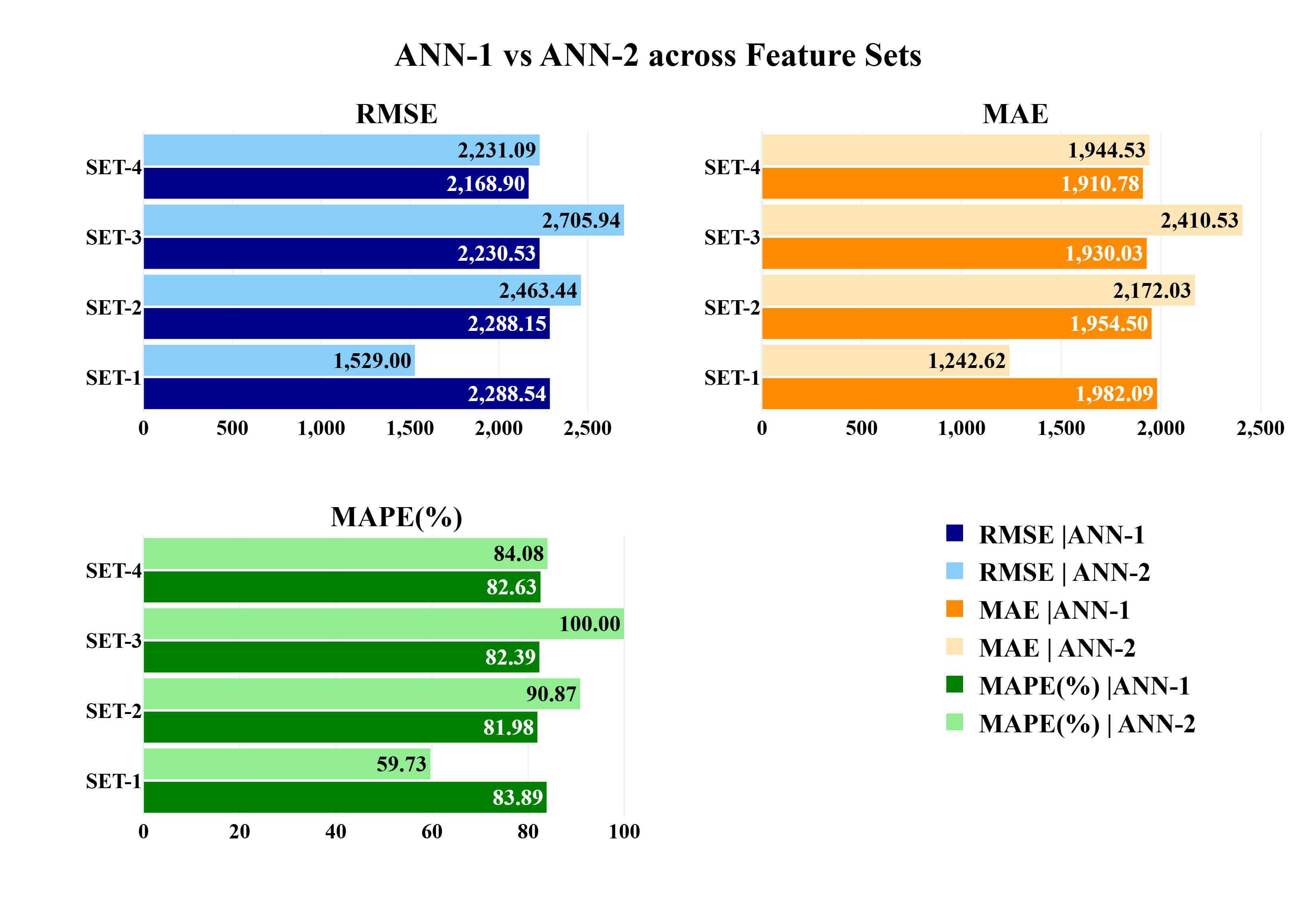}
	\caption{ANN results for the Barishal division, showing the error metrices (RMSE, MAE, MAPE) of the best performing model within each feature set.}	
	\label{annbe}
\end{figure}

ANN performance is division-specific, with different optimal configurations in Dhaka and Barishal. In Dhaka (Figure~\ref{annd}(a-b)), the lowest test error ANN-1 run is ANN-1 with SET-1 (temperature, rainy days, sunshine hours, humidity)(RMSE=2176.70, MAE=1282.00, MAPE=31.54$\%$) (Figure~\ref{annde}), indicating that dengue variability is well explained by the combined effects of thermal conditions and moisture availability, with frequent rainfall events for repeated breeding site refilling and sunshine hours modulating habitat persistence. The lowest error ANN-2 run is also SET-1 plus case-lag (RMSE=2322.90, MAE=1340.00, MAPE=28.53$\%$), suggesting that adding short-term persistence (1 month case lag) can improve proportional fit (lower MAPE) even when overall magnitude errors (RMSE, MAE) are slightly higher.
For Barishal (Figure~\ref{annb}(a-b)), the best performance is obtained by ANN-2 with SET-1 (temperature, rainy days, sunshine hours, humidity) in terms of RMSE (RMSE=1529.00, MAE=1242.62, MAPE=59.73$\%$) (Figure~\ref{annbe}), suggesting that short-term persistence (via the lagged dengue term) provides additional predictive power beyond climate alone, while the same core climate drivers still dominate seasonal suitability.
Among the climate-only ANN-1 runs, the best performance is obtained with SET-4 (temperature, rainfall, sunshine hours, humidity) (RMSE=2168.90), suggesting that rainfall amount (water availability for larval habitats) is informative for Barishal when case-lag is not used while sunshine hours and humidity modulate how long breeding habitats persist. However, its larger errors compared with the overall best ANN-2 result show that case persistence adds substantial predictive value for Barishal in the test window.
\subsection{XGBoost}
\begin{figure}[H]
	\centering 
	\subfloat[]{\includegraphics[width=4.2 in]{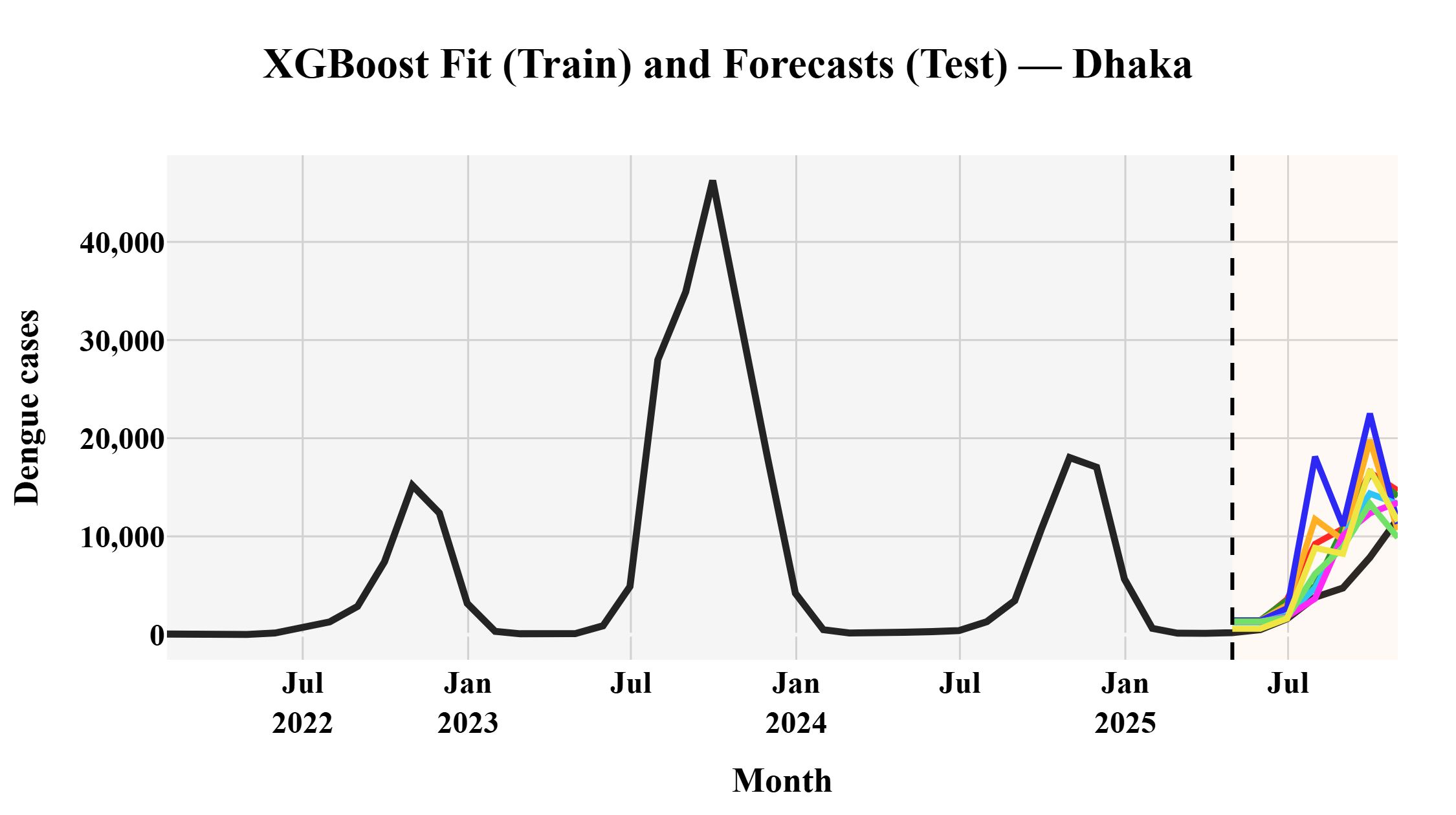}}\\
	\subfloat[]{\includegraphics[width=4.2 in]{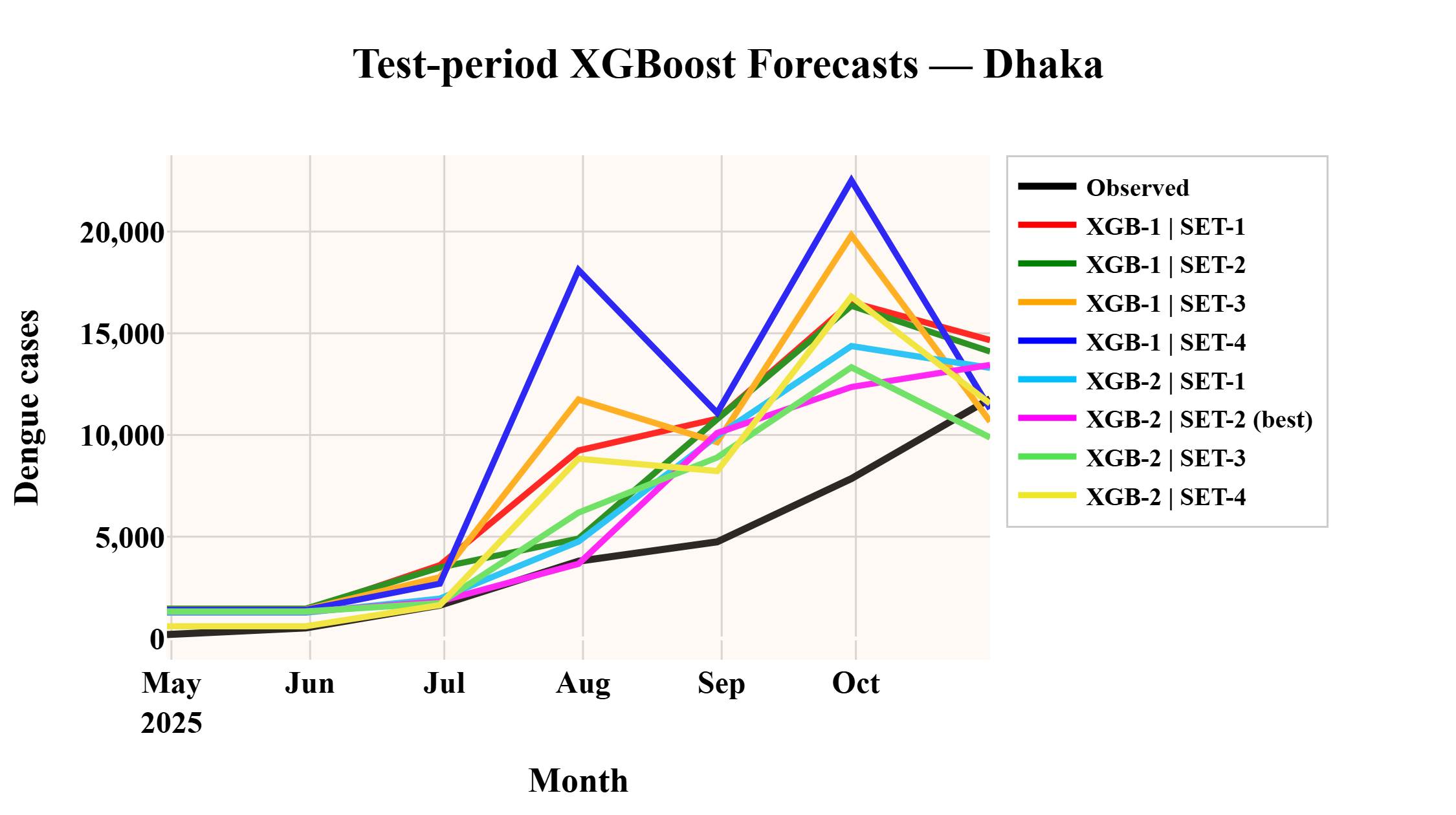}}
	
	\caption{XGBoost results for Dhaka division. (a) Full series observed cases with test-window forecasts from the best model in each feature set; (b) Test period observed vs forecasts (best per feature set).}	
	\label{xgbd}
\end{figure}

\begin{figure}[H]
	\centering 
	
	\includegraphics[width=6 in]{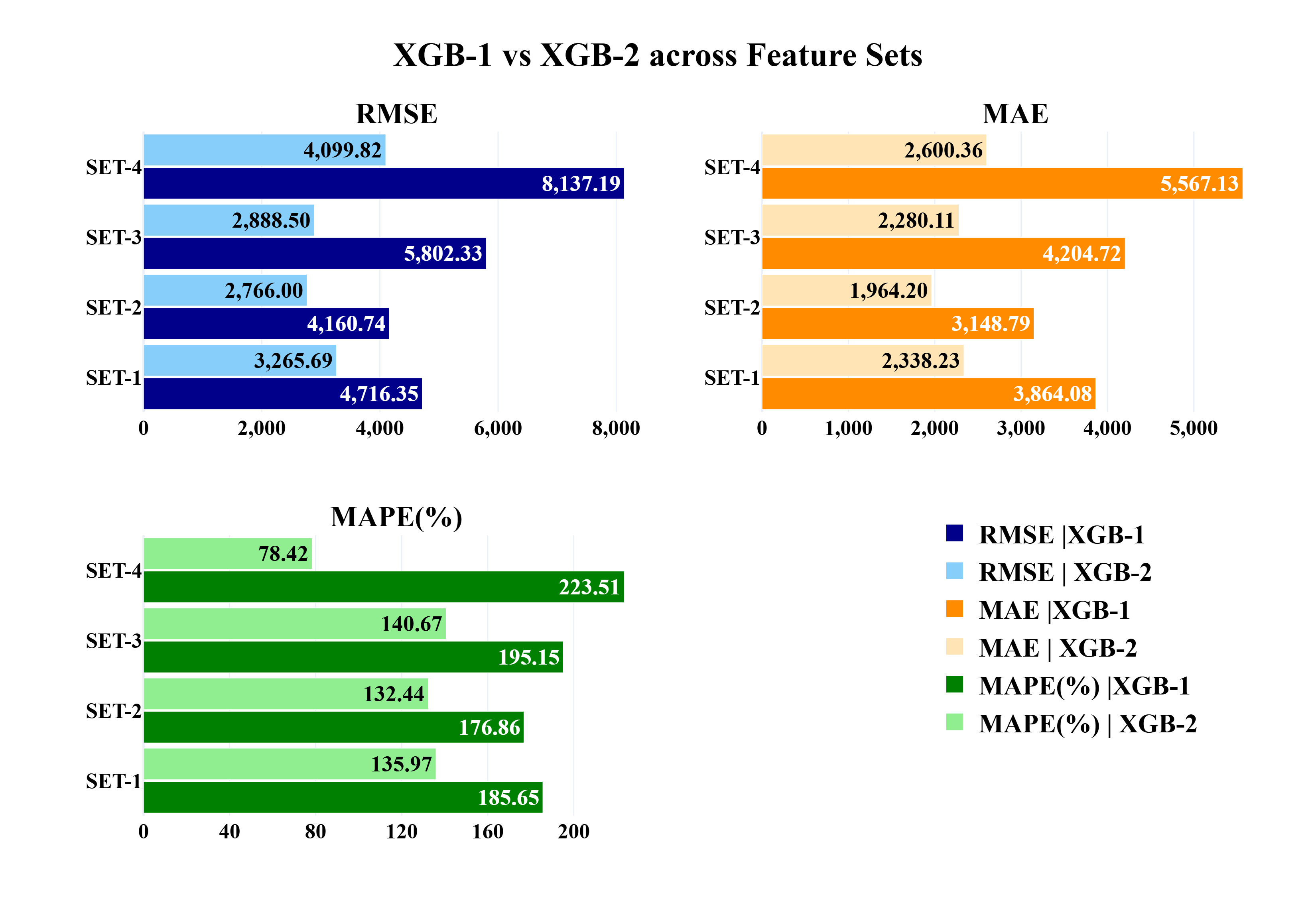}
	\caption{XGBoost results for the Dhaka division, showing the error metrices (RMSE, MAE, MAPE) of the best performing model within each feature set.}	
	\label{xgbde}
\end{figure}

\begin{figure}[H]
	\centering 
	\subfloat[]{\includegraphics[width=4.2 in]{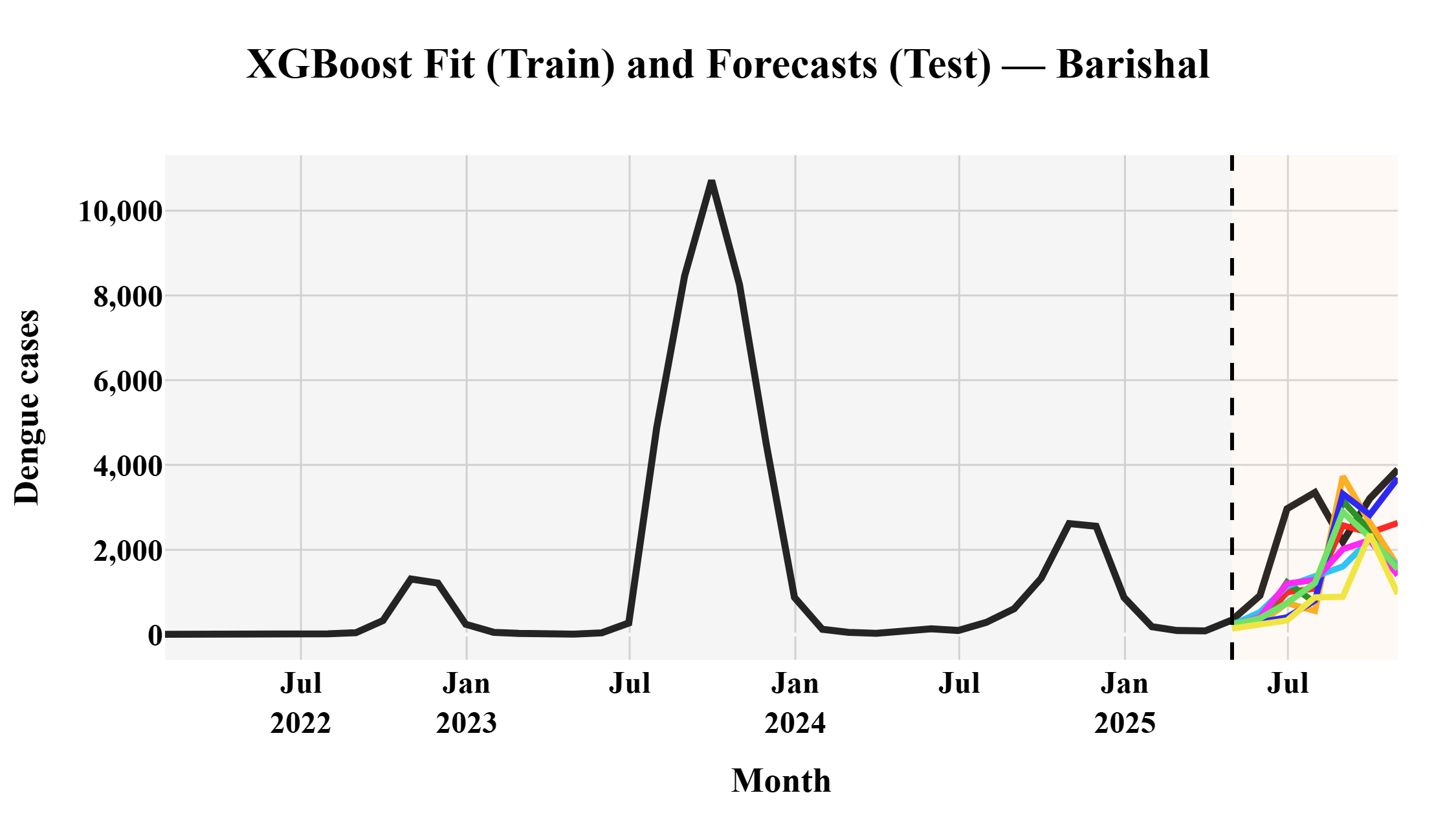}}\\
	\subfloat[]{\includegraphics[width=4.2 in]{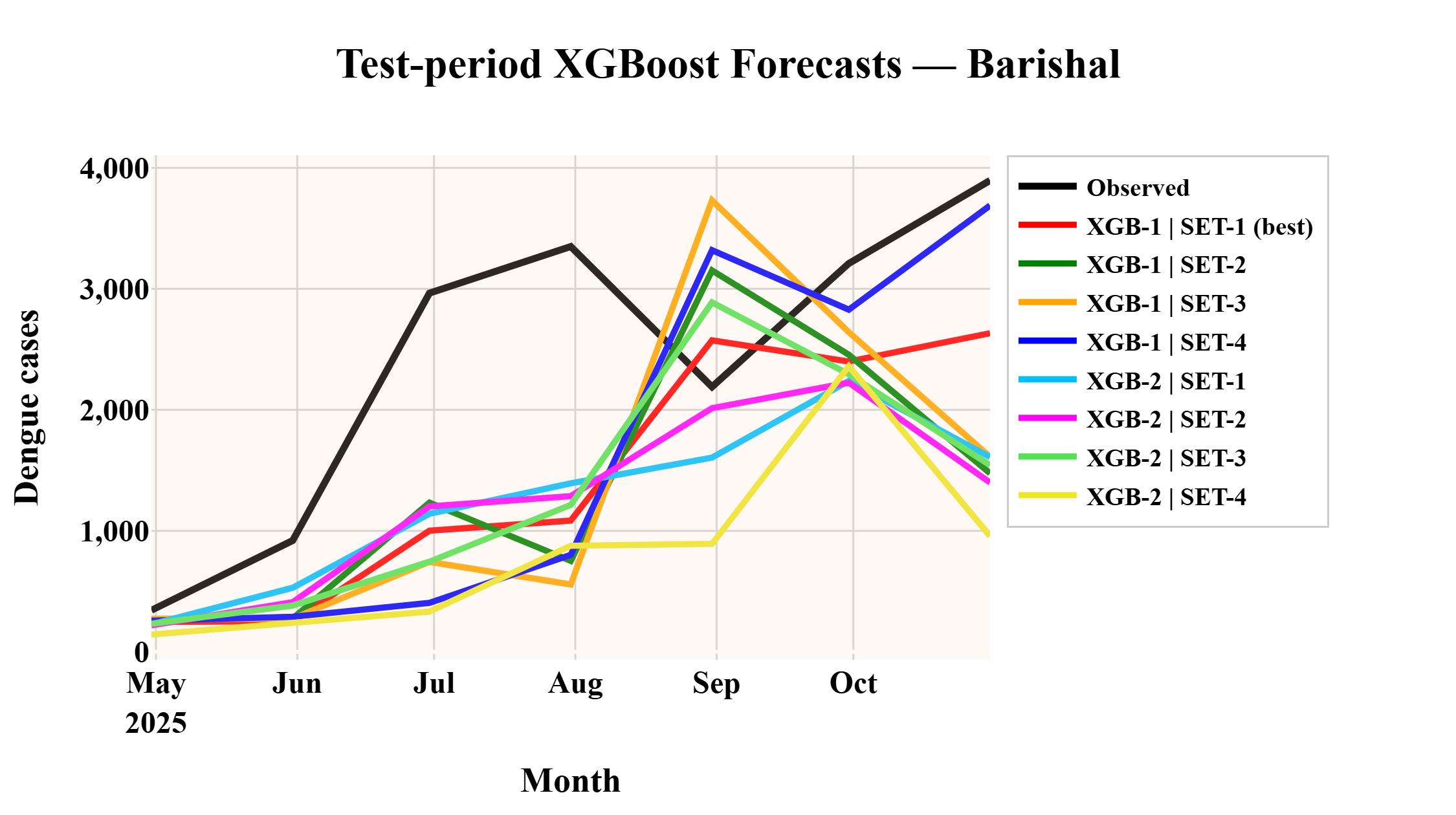}}
	
	\caption{XGBoost results for Barishal division. (a) Full series observed cases with test-window forecasts from the best model in each feature set; (b) Test period observed vs forecasts (best per feature set).}	
	\label{xgbb}
\end{figure}

\begin{figure}[H]
	\centering 
	
	\includegraphics[width=6 in]{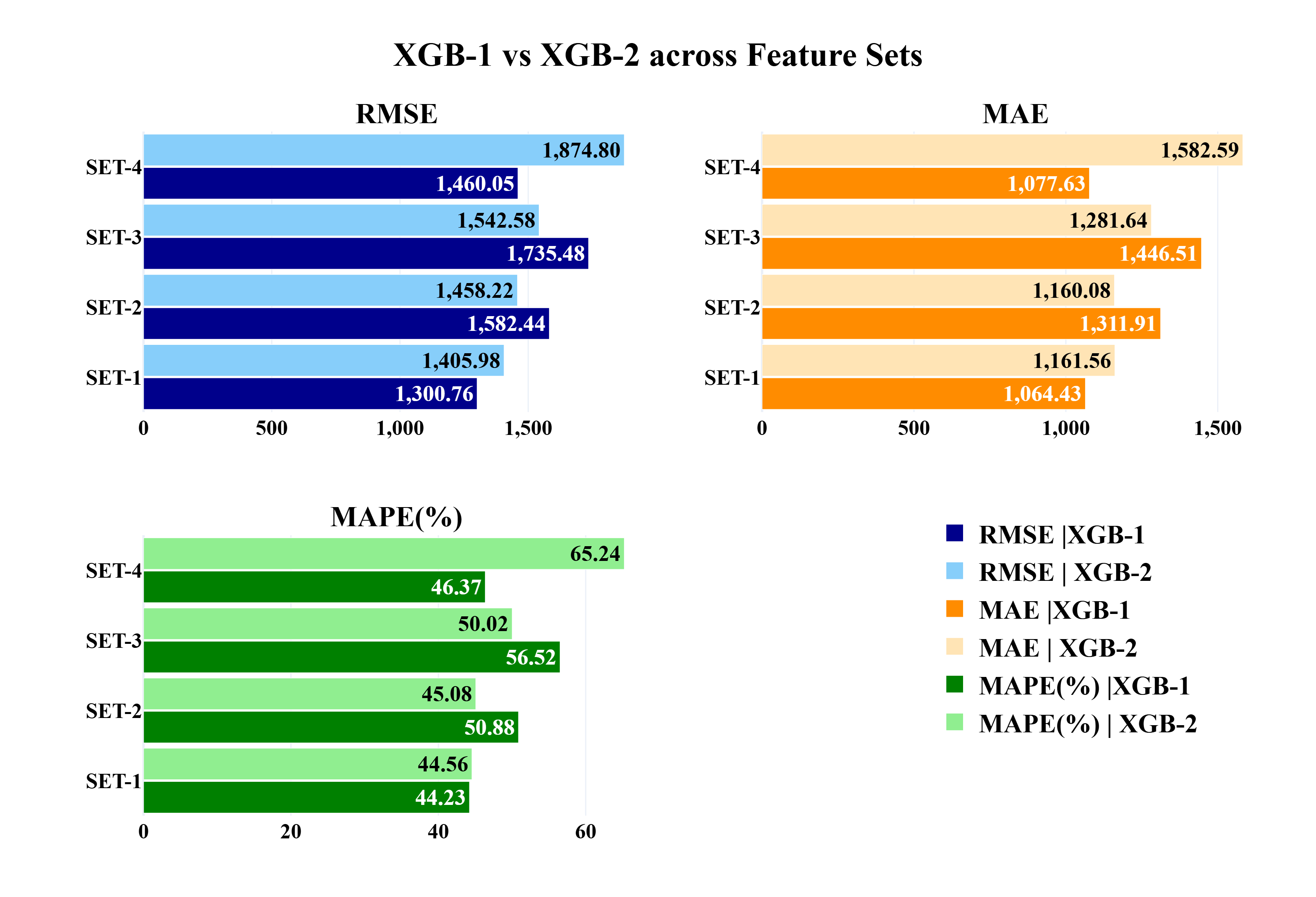}
	\caption{XGBoost results for the Barishal division, showing the error metrices (RMSE, MAE, MAPE) of the best performing model within each feature set.}	
	\label{xgbbe}
\end{figure}

Across divisions, the XGBoost results show a clear dissimilarity: Dhaka benefits more from adding case persistence (case-lag), whereas Barishal is explained well by climate-only features.\\
In Dhaka (Figure~\ref{xgbd}(a-b)), the best-performing XGBoost model is the case-lag (lag 1) version (XGB-2), with the lowest test error obtained by SET-2(Temp+Rainy Days+Sun Days+Humidity; RMSE=2766.0, MAE=1964.20, MAPE=132.44$\%$). The corresponding best climate-only option is XGB-1 with SET-2 (RMSE=4160.7, MAE=3148.79, MAPE=176.86$\%$) (Figure~\ref{xgbde}), so including the lagged dengue term noticeably improves short-term prediction. This suggests that, beyond climate forcing, Dhaka's dengue trajectory during the test window has stronger temporal momentum, while the climate signals in SET-2 point to roles for temperature and humidity (mosquito development and survival), rainy days (recurrent breeding-site replenishment), and sun days (proxy for higher solar exposure and faster drying of breeding sites). 
In Barishal (Figure $\ref{xgbb}$(a-b)), the pattern is different: the overall best model is climate-only XGB-1, specifically SET-1 (Temp+Rainy Days+Sunshine Hours+Humidity; RMSE=\allowbreak1300.76, MAE=1064.43, MAPE=44.23~\%). The best case-lag counterpart, XGB-2 with SET-1, is slightly worse (RMSE=1405.99, MAE=1161.56, MAPE=44.56$\%$) (Figure~\ref{xgbbe}), implying that case persistence adds limited value here relative to the strength of climate forcing alone. 
The best Barishal feature set emphasizes temperature and humidity together with rainy days and sunshine hours, where sunshine can act as a proxy for drying and evaporation and thus modulate how effectively rainfall translates into sustained larval habitats.
\subsection{Division-wise Best Model Selection}
\begin{figure}[H]
	\centering 
	\subfloat[]{\includegraphics[width=4.5 in]{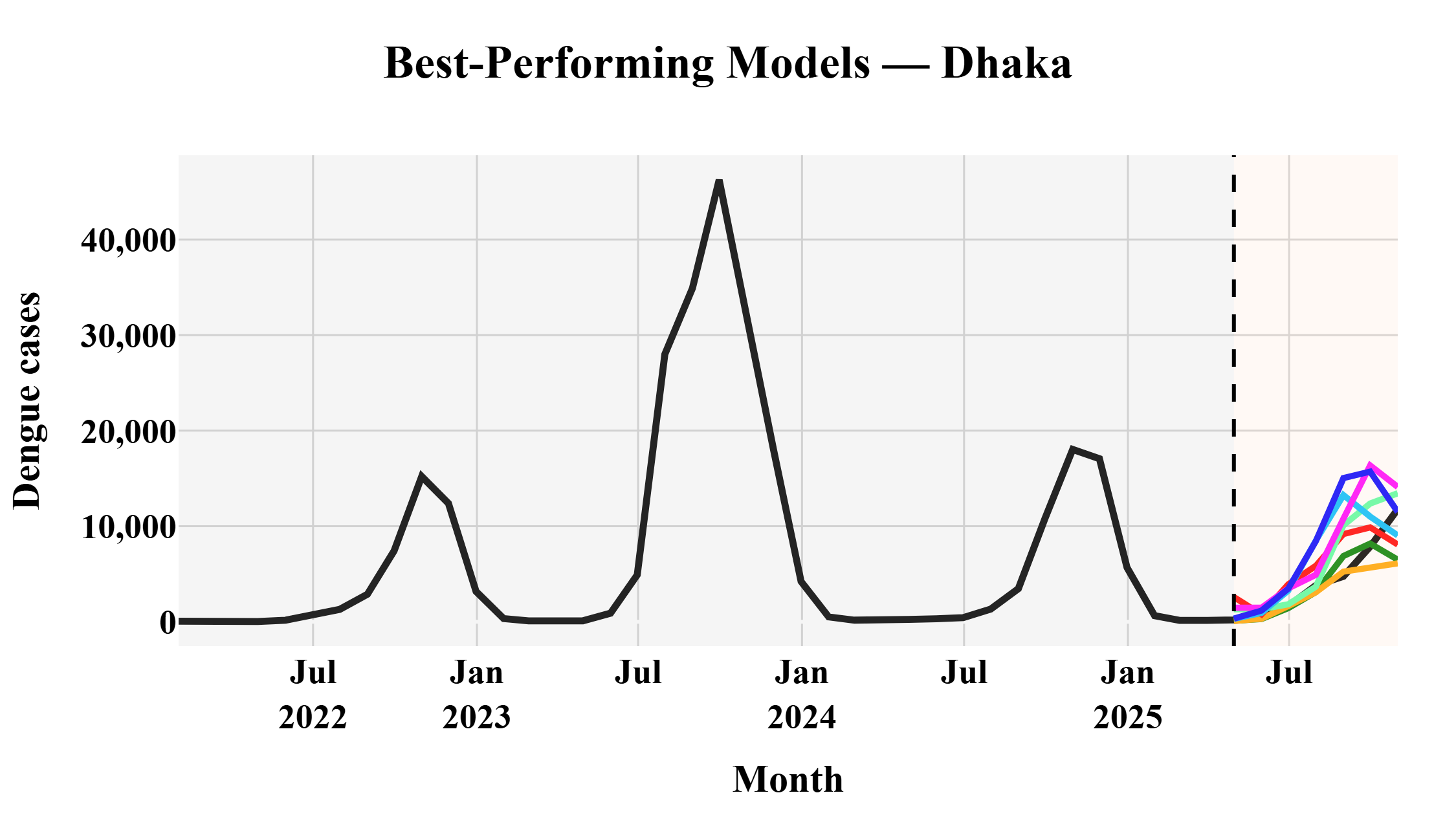}}\\
	\subfloat[]{\includegraphics[width=4.5 in]{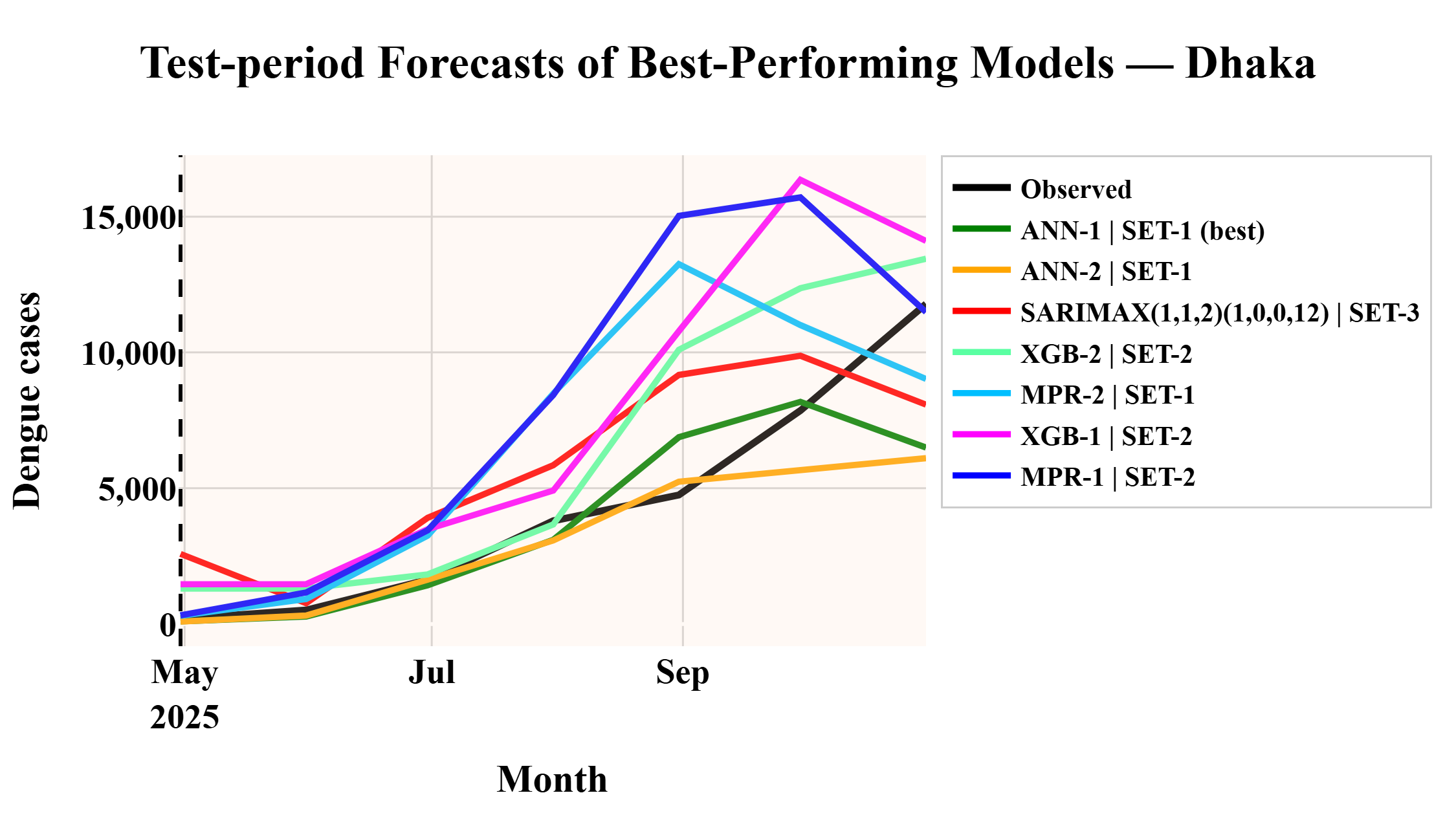}}
	\caption{Best-performing models by feature set for Dhaka division (a) Full series with the test window highlighted and best-model forecasts overlaid; (b) Test-period zoom showing observed cases versus the corresponding best model forecasts.}
	\label{bestd}
	\end{figure}

\begin{figure}[H]
	\centering 
	\subfloat[]{\includegraphics[width=4.5 in]{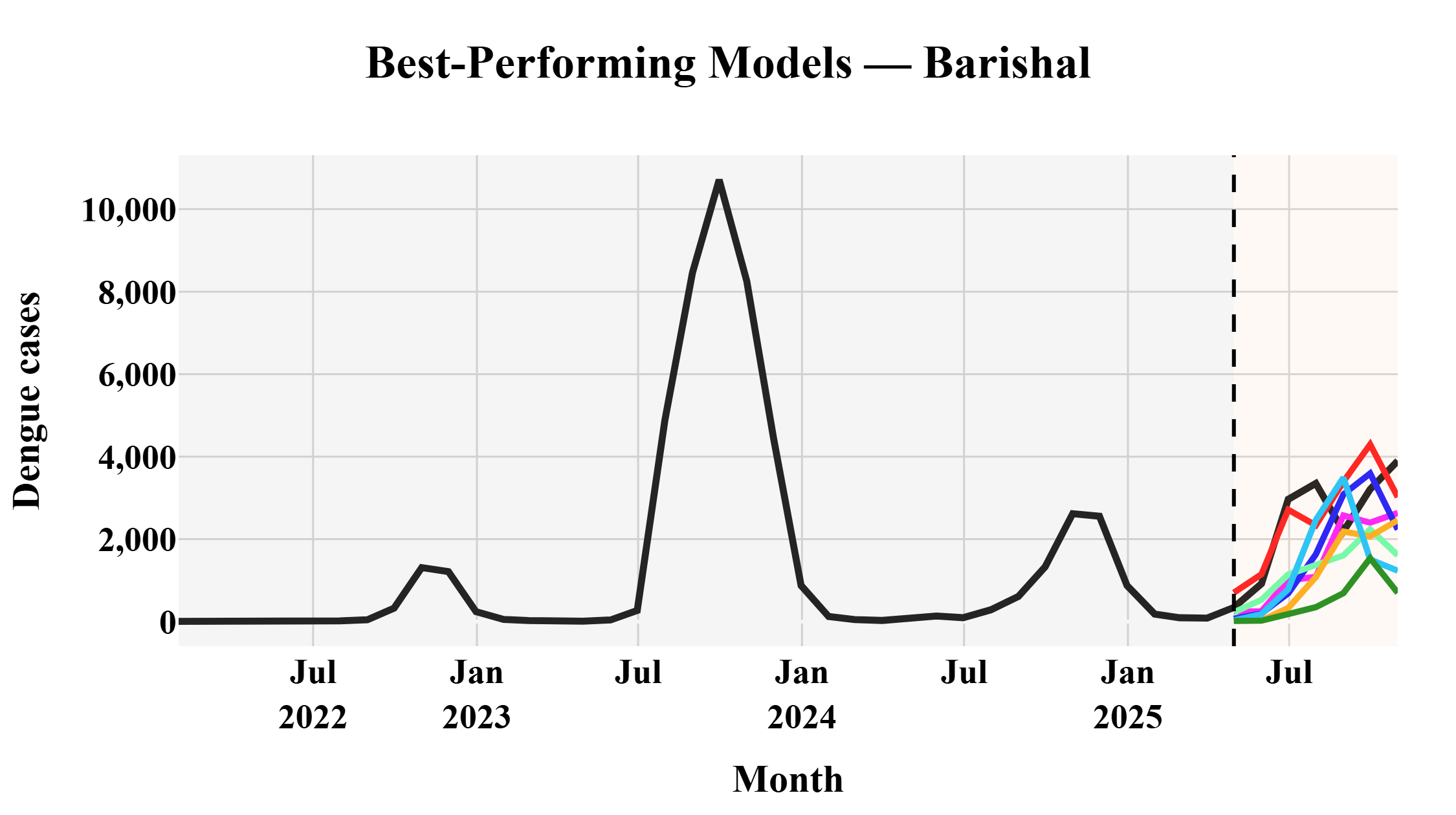}}\\
	\subfloat[]{\includegraphics[width=4.5 in]{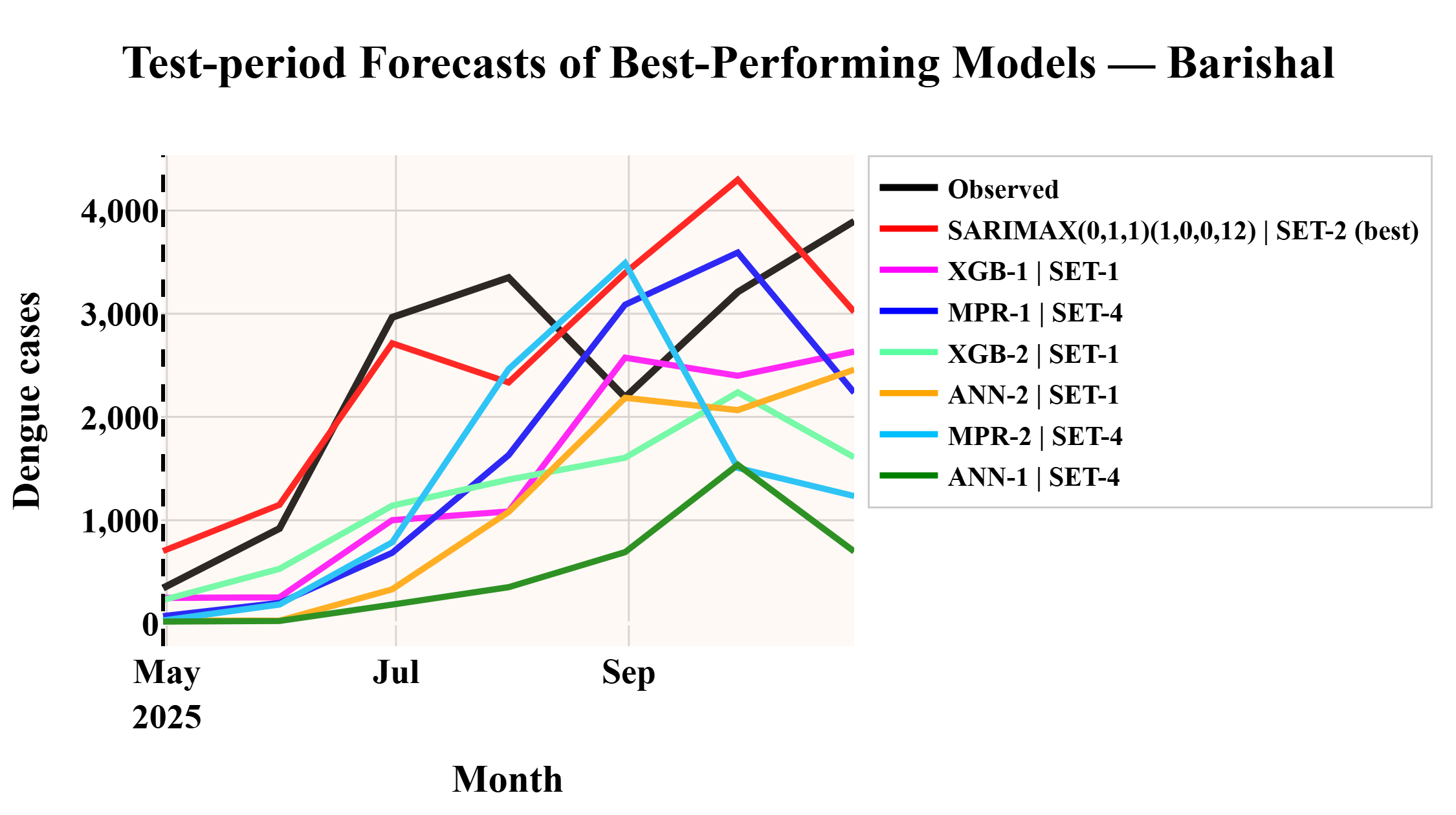}}
	\caption{Best-performing models by feature set for Barishal division (a) Full series with the test window highlighted and best-model forecasts overlaid; (b) Test-period zoom showing observed cases versus the corresponding best model forecasts.}	
	\label{bestb}
\end{figure}

Here we only compare the best selected model from each model class (after testing all four feature sets earlier) over the test window and then identify the best among these best models for each division.\\
For Dhaka (Figure~\ref{bestd}(a-b)), the best overall performance among the selected models is obtained by ANN-1 with SET-1 (temperature, rainy days, sunshine hours, humidity). This indicates that dengue variation is most consistently explained by the combined effects of temperature and humidity, together with rainy days (repeated replenishment of breeding sites) and sunshine hours (drying conditions that influence how long standing water persists).
For Barishal (Figure~\ref{bestb}(a-b)), the best overall model among the selected models is SARIMAX with SET-2 (temperature, rainy days, sun days, humidity) with annual seasonality. This indicates that, in Barishal, dengue patterns are explained best by a seasonal SARIMAX structure with lagged climate predictors, whereas the case-lag variants do not provide the lowest error in this comparison.
From a public-health perspective, the key information to focus on is a consistent early warning pattern: frequent rainy days followed by sustained temperature and humidity conditions that support mosquito breeding and survival. In settings where vaccination is not available, prevention depends mainly on vector control and early action, so control should be strengthened before case counts rise by increasing source reduction in households and public spaces, maintaining drainage and removing water-holding containers, applying larval control where water cannot be eliminated, and intensifying vector surveillance with rapid targeted response in the highest-risk areas during these periods.
\section{Discussion and Conclusions}\label{Discussion}
Dengue transmission in Bangladesh is strongly seasonal and aligns closely with the monsoon--post-monsoon window, when temperature-humidity conditions and water availability favor \textit{Aedes} breeding. In both Dhaka and Barishal, cases remain low from January to May, rise from June, and typically peak around September--October. This timing is consistent with the synchronized climate cycle observed in both divisions: monsoon months bring high humidity, frequent rainy days and high rainfall, while sunshine (sun hours and sun days) drops sharply. The largest outbreak signal appears in September 2023, following the core monsoon period, supporting the expected lag between intensified vector breeding conditions and reported dengue incidence.\\
STL decomposition confirms clear annual seasonality in both divisions, with Dhaka showing a larger seasonal amplitude than Barishal. Dhaka's trend rises modestly up to about 2023--early 2024 and then declines, whereas Barishal's trend increases steadily. Large 2023 residual spikes indicate outbreak shocks beyond seasonality, and ACF--PACF results show strong short-term persistence, supporting inclusion of lagged dengue incidence as an additional predictor in selected model configurations. Lagged correlation analysis provides a coherent picture of delayed climate forcing across both divisions: rainfall-related variables (rainy days and rainfall) show the strongest positive association at about a two month lag, humidity peaks positively at about a one- month lag, sunshine-related measures (sun days and sunshine hours) are most strongly negative around a two-month lag, and temperature exhibits weaker but positive association at longer lags (three to four months). Building on this, a key contribution of the study is the systematic comparison of four climate feature sets that alternate wetness and sunshine representations while keeping temperature and humidity fixed, revealing that the most informative formulation is division-specific, because rainy days vs. rainfall and sun days vs. sun hours capture different aspects of habitat refilling, water volume, cloudiness, and drying potential that shape mosquito suitability at the monthly scale.\\
Model comparisons show that forecast skill depends on both model class and climate feature representation, and the optimal choice is division-specific. In Dhaka, ANN-1 with SET-1 (temperature, rainy days, sunshine hours, humidity) performs best, suggesting that nonlinear links between dengue and combined wetness-sunshine conditions are important, and that rainfall frequency and sunshine duration better capture habitat persistence and drying than rainfall totals or sun-day counts. In Barishal, SARIMAX with SET-2 (temperature, rainy days, sun days, humidity) performs best, indicating that seasonal structure with well-chosen exogenous regressors captures dynamics more effectively, and that sun-day counts may better reflect sustained cloudy conditions and habitat persistence. Comparing predictor structures further shows that adding a case-lag is more beneficial in Dhaka (stronger short-term persistence), whereas climate-only formulations are often sufficient or better in Barishal (stronger direct climate forcing), implying that dengue predictability reflects both climate forcing and intrinsic persistence, with their relative roles varying by division.\\

From a public-health perspective, the stable lag structure--about two months for wetness and sunshine measures and about one month for humidity--provides a practical early-warning window, where monitoring monsoon rainfall and reduced sunshine can guide timely vector control, risk communication, and preparedness ahead of the usual September--October peak. The fact that the best feature sets and models differ between Dhaka and Barishal also indicates that forecasting systems should be locally calibrated rather than using a single national specification. In addition, while current dengue mitigation in Bangladesh mainly relies on vector control and personal protection with clinical management, vaccination may serve as a valuable complementary option if safe, effective, and operationally feasible programs become available \cite{Khondaker2025acta}; in that case, improved forecasts could help identify where and when preventive efforts should be intensified. More generally, forecast-based alerts can be operationalized by mapping predicted risk to time-varying intervention intensities and resource-allocation schedules, consistent with optimal-control approaches developed across infectious-disease settings \cite{Khondaker2025,Khondaker2022abr}. The study shows that effective dengue forecasting in Bangladesh depends not only on model choice but also on how climate predictors and lags are constructed, supporting tailored early-warning strategies and motivating future extensions that incorporate additional drivers and longer histories to improve reliability during major outbreak years.

\section*{Acknowledgments}
%The authors’ acknowledged to the anonymous reviewers for their suggestions which significantly improved the quality of the manuscript. 
The research was partially supported by the University Grants Commission (UGC), %the Ministry of Science and Technology,
 Bangladesh. %Bose Center for Advanced Study and Research in Natural Sciences,  for supporting students to complete their graduation smoothly.

%=========================Conflict of interest======================================
\section*{Conflict of interest}
The authors declare no conflict of interest. 

\section*{Data Statement}
The datasets generated and/or analysed during the current study are available in the Faizunnesa Khondaker repository, \url{https://github.com/Faizunnesa-Khondaker/Dengue.git}.
 %https://github.com/Faizunnesa-Khondaker/Dengue.git
 
 \section*{Ethical approval}
 No consent is required to publish this manuscript.
 
 \section*{Author Contributions (CRediT)}
 Faizunnesa Khondaker: Conceptualization, Data curation, Formal analysis, Methodology, Software, Validation, Visualization, Writing–original draft.\\
 Md. Kamrujjaman:  Methodology, Investigation, Validation, Software, Supervision, Writing–review and editing.\\
 All authors have read and agreed to the published version of the manuscript.

\end{document}